\batchmode
\documentclass[11pt]{article}

\usepackage{epsfig}
\usepackage{graphicx}
\usepackage{color}

\newtheorem{theorem}{Theorem}
\newtheorem{lemma}{Lemma}
\newtheorem{corollary}{Corollary}

\newtheorem{proposition}{Proposition}
\newtheorem{definition}{Definition}

\newcommand{\be}{\begin{equation}}
\newcommand{\ee}{\end{equation}}
\newcommand{\bea}{\begin{eqnarray}}
\newcommand{\eea}{\end{eqnarray}}
\newcommand{\beas}{\begin{eqnarray*}}
\newcommand{\eeas}{\end{eqnarray*}}
\newcommand{\ba}{\begin{array}}
\newcommand{\ea}{\end{array}}

\newcommand{\real}{\mbox{$\mathrm{I\!R}$}}

\newcommand{\eps}{\ensuremath{\epsilon}}

\newcommand{\bfD}{\ensuremath{\mathbf{D}}}

\def\XXint#1#2#3{{\setbox0=\hbox{$#1{#2#3}{\int}$}
     \vcenter{\hbox{$#2#3$}}\kern-.5\wd0}}

\newcommand{\mcG}{\ensuremath{\mathcal{G}}}

\newcommand{\mcL}{\ensuremath{\mathcal{L}}}

\newcommand{\mcP}{\ensuremath{\mathcal{P}}}

\def\qed{\hbox{\vrule width 6pt height 6pt depth 0pt}}

\usepackage{amsmath}
\usepackage{amssymb}

\title{Regularity of the Solution to 1-D Fractional Order Diffusion Equations} 
\author{V.~J.~Ervin\thanks{Department of Mathematical Sciences,
	  Clemson University, Clemson, South Carolina 29634-0975.
	  email: {\tt vjervin@clemson.edu}. Partially support by 
   CONICYT through FONDECYT project 1150056.}
	\and 
	N.~Heuer\thanks{Facultad
	  de Matem\'aticas, Pontificia Universidad Cat\'olica de Chile,
	  Avenida Vicu\~na Mackenna 4860, Macul, Santiago, Chile.
	  email: {\tt nheuer@mat.puc.cl}. Partially support by 
   CONICYT through FONDECYT projects 1150056, 
   and Anillo ACT1118 (ANANUM).}
	\and
	J.~P.~Roop\thanks{Department of Mathematics, 
	North Carolina A \& T State University, Greensboro,
      North Carolina 27411.
	  email: {\tt jproop@ncat.edu}.}
	}

\date{07/28/16}

\begin{document}
\maketitle

\begin{abstract}
In this article we investigate the solution of the steady-state fractional diffusion equation on a bounded 
domain in $\real^{1}$. From an analysis of the underlying model problem, we postulate that the fractional
diffusion operator in the modeling equations is neither the Riemann-Liouville nor the Caputo fractional
differential operators. We then find a closed form expression for the kernel of the fractional diffusion
operator which, in most cases, determines the regularity of the solution. Next we establish that the
Jacobi polynomials are pseudo eigenfunctions for the fractional diffusion operator. A spectral type
approximation method for the solution of the steady-state fractional diffusion equation is then 
proposed and studied. 
\end{abstract}



\setcounter{equation}{0}
\setcounter{figure}{0}
\setcounter{table}{0}
\setcounter{theorem}{0}
\setcounter{lemma}{0}
\setcounter{corollary}{0}
\section{Introduction}
 \label{sec_intro}
The history of the fractional derivative is almost as long as the history of the integer order derivative.
Notably from a 1695 letter of Leibniz to L'H\^{o}pital referring to the question of fractional order derivatives,
``Il y a de l'apparence qu'on tirera un jour des
consequences bien utiles de ces paradoxes, car il n'y a gueres de paradoxes sans utilit\'{e}, '' 
which translates to ``It will lead to a paradox, from which one day useful consequences will be drawn.''
 In recent years the fractional 
derivative has received increased attention in modeling a variety of physical phenomena. Most often cited
are applications in contaminant transport in ground water flow \cite{ben001},
viscoelasticity \cite{mai971}, turbulent flow \cite{mai971, shl871}, and chaotic dynamics \cite{zas931}.
As interest in the fractional derivative has increased so has approximation methods
to solve such equations. Generally speaking (for the 1-D case), approximation methods which exist for
integer order differential equations have been successfully adapted to the fractional order case. Specifically,
to mention a few (a complete list is beyond the focus of this article), 
finite difference methods \cite{cui091, liu041, mee041, tad071, wan121},
finite element methods \cite{erv061, jin151, liu111, wan131}, 
discontinuous Galerkin methods \cite{xu141}, mixed methods \cite{che161},
spectral methods \cite{che162, li121, wan151, zay151}, enriched subspace methods \cite{jin161}. 
To date most of the approximation
schemes have focused on the 1-sided fractional diffusion equation
\be
 \mcL_{1}^{\alpha} u(x) \ := \ - \bfD^{\alpha} u(x) \ = \ f(x) , \ x \in (0 , 1), \ \ u(0) = u(1) = 0 \, , 
  \label{reg1}
\ee
for $1 < \alpha < 2$. ( A formal definition of $\bfD^{\alpha} u(x)$ is given in the following section.)

Another interesting historical fact, a point of particular interest in this article, is the definition of the
fractional derivative. Or more precisely stated, \textit{definitions} of the fractional derivative.
There has been a number of definitions of the fractional derivatives studied. Most relevant to our discussion
are the Riemann-Liouville fractional derivative and the Caputo fractional derivative. We refer the reader 
to the monographs \cite{bal121, kil061, pod991, sam931} for a detailed discussion of various fractional 
derivatives. Also, of particular note is the recent approach to modeling nonlocal diffusion problems using a
linear integral operator introduced by Du, Gunzburger, Lehoucq et al. (see \cite{du121}).

Motivated by our
interest in physical applications, in the following section we present the Riemann-Liouville and Caputo 
fractional derivatives on a finite interval, which for the sake of specificity we take to be $I \, := \, (0 , 1)$.
(In the case where a function and its (integer) derivatives vanish at the endpoint of the interval the
Riemann-Liouville and Caputo fractional derivatives agree.)
 
The motivation for this article was to investigate the regularity of the solution to the two-sided
fractional diffusion equation
\be
 \mcL_{r}^{\alpha} u(x) \ := \ - \left( r \bfD^{\alpha} u(x) \ + \ (1 - r) \bfD^{\alpha *} u(x) \right) 
 \ = \ f(x) , \ x \in (0 , 1), \ \ u(0) = u(1) = 0 \, , 
  \label{reg2}
\ee
for $1 < \alpha < 2$, and $0 < r < 1$, which we think is a more physical model of diffusion than \eqref{reg1}.
(In \eqref{reg2} diffusion occurs to both the left and right of any point in the domain.) A variational formulation
of the solution to \eqref{reg2} was studied in \cite{erv061}, together with a finite element error analysis.
The error analysis was based on assumptions on the regularity of the true solution $u$, which has been
pointed out by a number of other authors, is not justified for a general right hand side function $f$.
In \cite{jin151} Jin et al. presented a very nice analysis and discussion of the regularity of the solution
to \eqref{reg1} for $\bfD^{\alpha}$ interpreted as the Riemann-Liouville fractional derivative
and as the Caputo fractional derivative.
In general, the solution of \eqref{reg1} has a singularity in the derivative at $x = 0$. Very helpful in studying
the regularity of the solution to \eqref{reg1} is the existence of an explicit inverse to $\mcL_{1}^{\alpha}$
which satisfies $\left( \mcL_{1}^{\alpha} \right)^{-1} f (0) \, = \, 0$. We do not have an explicit inverse
for $\mcL_{r}^{\alpha}$. Subsequently we have to think more generally about the operator $\mcL_{r}^{\alpha}$,
and in particular the definition of $\bfD^{\alpha}$ in the context of diffusion problems.

Following the introduction of notation in Section \ref{sec_not},  in Section \ref{sec_mod} we present a discussion 
on the modeling of the fractional diffusion equation. We subsequently conclude that in the context of a diffusion
operator the appropriate interpretation of the fractional derivative is neither the Riemann-Liouville definition nor the
Caputo definition. Rather, for $1 < \alpha < 2$,
\be
   \bfD^{\alpha} u(x) \ := \ D \, \bfD^{-(2 - \alpha)} \, D u(x) \, .
  \label{reg3}
\ee  

The kernel of the operator $\mcL_{r}^{\alpha}$, $ker(\mcL_{r}^{\alpha})$
plays a key role in determining the regularity of the solution of \eqref{reg2}.
Thus the definition of $\bfD^{\alpha}$
is central in determining the regularity of the solution to \eqref{reg2}. In Section \ref{sec_reg}
we discuss the regularity of the solution to \eqref{reg2}, using the definition of $\bfD^{\alpha}$
given in \eqref{reg3}. Somewhat of a surprise is that the regularity of the solution depends
upon $r$. 
In order to numerically illustrate the regularity of the solution to \eqref{reg2} 
in Section \ref{secFEMcvg} we present Finite Element Method (FEM) computations.
The experimental rates of convergence of the FEM approximations are consistent
with the regularity of the solution obtained in Section \ref{sec_reg}. 

In Section \ref{secSpec} we establish that Jacobi polynomials are pseudo eigenfunctions for the 
fractional diffusion operator. Specifically (see Lemmas \ref{lmaGnO} and \ref{lmaGngO}) we show that
\[
    \mcL_{r}^{\alpha} \omega(x) \, \mcG_{n}(x) \ = \ \lambda_{n} \, \mcG_{n}^{*}(x) \, ,
\]
where  $\mcG_{n}(x)$ and  $\mcG_{n}^{*}(x)$ are Jacobi polynomials,
$\omega(x)$ is the Jacobi weight, and $\lambda_{n}$ the pseudo eigenvalue. Using this property we propose
and study a spectral type approximation method for the solution of steady-state fractional diffusion equations.
Two numerical examples are given to illustrate the performance of the method.

\setcounter{equation}{0}
\setcounter{figure}{0}
\setcounter{table}{0}
\setcounter{theorem}{0}
\setcounter{lemma}{0}
\setcounter{corollary}{0}
\section{Notation and Properties}
\label{sec_not}
For $u$ a function defined on $(a , b)$, and $\sigma > 0$,
we have that the left and right fractional integral operators are defined as: \\
\underline{Left Fractional Integral Operator}:
 $ \mbox{}_{a}D_{x}^{-\sigma}u(x) \, := \, \frac{1}{\Gamma(\sigma)} \int_{a}^{x} (x - s)^{\sigma - 1} \, u(s) \, ds \, . $ 
 \vspace{0.5em} \\
\underline{Right Fractional Integral Operator}:
 $ \mbox{}_{x}D_{b}^{-\sigma}u(x)  \, := \, \frac{1}{\Gamma(\sigma)} \int_{x}^{b} (s - x)^{\sigma - 1} \, u(s) \, ds \, . $ 
 \vspace{0.5em} \\
Then, for $\mu > 0$, $n$ the smallest integer greater than $\mu$ $(n-1 \le \mu < n)$,  
$\sigma = n - \mu$, and $D$ the derivative operator,
the left and right  Riemann-Liouville fractional differential operators are defined as:  \\
\underline{Left Riemann-Liouville Fractional Differential Operator of order $\mu$}: \\
 \mbox{} \hfill $ \mbox{}_{a}^{RL}D_{x}^{\mu} u(x) \, := \, D^{n} \mbox{}_{a}D_{x}^{-\sigma}u(x) \ = \ 
 \frac{1}{\Gamma(\sigma)} \frac{d^{n}}{dx^{n}} \int_{a}^{x} (x - s)^{\sigma - 1} \, u(s) \, ds \, . $ 
 \vspace{0.5em} \\
\underline{Right Riemann-Liouville Fractional Differential Operator of order $\mu$}: \\
 \mbox{} \hfill  $ \mbox{}_{x}^{RL}D_{b}^{\mu} \, u(x) := \, (- D)^{n} \mbox{}_{x}D_{b}^{-\sigma}u(x) \ = \ 
 \frac{(-1)^{n}}{\Gamma(\sigma)} \frac{d^{n}}{dx^{n}} \int_{x}^{b} (s - x)^{\sigma - 1} \, u(s) \, ds \, . $ 
 
The Riemann-Liouville and Caputo fractional differential operators differ in the location of the derivative operator. \\
\underline{Left Caputo Fractional Differential Operator of order $\mu$}: \\
 \mbox{} \hfill $ \mbox{}_{a}^{C}D_{x}^{\mu} u(x) \, := \,  \mbox{}_{a}D_{x}^{-\sigma}\, D^{n} u(x) \ = \ 
 \frac{1}{\Gamma(\sigma)} \int_{a}^{x} (x - s)^{\sigma - 1} \, \frac{d^{n}}{ds^{n}} u(s) \, ds \, . $ 
 \vspace{0.5em} \\
\underline{Right Caputo Fractional Differential Operator of order $\mu$}: \\
 \mbox{} \hfill  $ \mbox{}_{x}^{C}D_{b}^{\mu} \, u(x) := \, (-1)^{n} \mbox{}_{x}D_{b}^{-\sigma}\, D^{n} u(x) \ = \ 
 \frac{(-1)^{n}}{\Gamma(\sigma)} \int_{x}^{b} (s - x)^{\sigma - 1} \, \frac{d^{n}}{ds^{n}} u(s) \, ds \, . $ 

As our interest is in the solution of fractional diffusion equations on a bounded, connected subinterval of $\real$, 
without loss
of generality we restrict our attention to the unit interval $(0 , 1)$.

For $s \ge 0$ let $H^{s}(0 , 1)$ denote the Sobolev space of order $s$ on the interval $(0 , 1)$, and
$\tilde{H}^{s}(0 , 1)$ the set of functions in $H^{s}(0 , 1)$ whose extension by $0$ are in 
$H^{s}(\real)$. Equivalently, for $u$ defined on $(0 , 1)$ and $\tilde{u}$ its extension by zero, 
$\tilde{H}^{s}(0 , 1)$ is the closure of $C_{0}^{\infty}(0 , 1)$ with respect to the norm 
$\| u \|_{\tilde{H}^{s}(0 , 1)} := \| \tilde{u} \|_{H^{s}(\real)}$.  With respect to
$L^{2}$ duality, for $s \ge 0$ we let $H^{-s}(0 , 1) := \left(\tilde{H}^{s}(0 , 1)\right)^{\prime}$, the dual space of 
$\tilde{H}^{s}(0 , 1)$.

Useful below in establishing results about the kernel of the fractional diffusion operator is the 
hypergeometric function \cite{leb721, sze751}.

\begin{definition} \label{defHGF}
The Gaussian three-parameter hypergeometric function $\mbox{}_{2}F_{1}$ is defined by an integral 
and series as follows:
\be
\mbox{}_{2}F_{1}(a, \, b; \, c; \, x) \ = \ 
\frac{\Gamma(c)}{\Gamma(b) \, \Gamma(c - b)} \int_{0}^{1} z^{b - 1} (1 - z)^{c - b - 1} (1 \, - \, z x)^{-a} \, dz
\ = \ \sum_{n = 0}^{\infty} \frac{ ( a )_{n} \, ( b )_{n} \, x^{n}}{( c )_{n} \, n!} \, ,
\label{eww3}
\ee
with convergence only if $Re(c) > Re(b) > 0$.
\end{definition}
In \eqref{eww3} $(q)_{n}$ denotes the (rising) Pochhammer symbol.

\begin{proposition} (Interchange property) \label{HGFp1}
For $Re(c) > Re(b) > 0$, and $Re(c) > Re(a) > 0$, we have that
\be
 \mbox{}_{2}F_{1}(a, \, b; \, c; \, x) \ = \ \mbox{}_{2}F_{1}(b, \, a; \, c; \, x) \, .
\label{eww4}
\ee
\end{proposition}
%


For ease of notation, we use 
\[
   \bfD^{-\sigma} \, := \,  \mbox{}_{0}D_{x}^{-\sigma} \, ,  \ \mbox{ and } \ 
    \bfD^{-\sigma *} \, := \,  \mbox{}_{x}D_{1}^{-\sigma} \, .
\]

\setcounter{equation}{0}
\setcounter{figure}{0}
\setcounter{table}{0}
\setcounter{theorem}{0}
\setcounter{lemma}{0}
\setcounter{corollary}{0}
\section{Interpretation of the Fractional Derivative}
 \label{sec_mod}
In this section we discuss the interpretation of the fractional derivative for modeling diffusion phenomena.

With \eqref{reg1} and \eqref{reg2} interpreted as the steady-state equation for a time dependent 
diffusion equation, we begin with a review of the derivation of the 1-D heat equation. 
(See \cite{goc111} for a more complete derivation.)

\begin{figure}[!ht]
\begin{center}
  \includegraphics[width=0.75\linewidth]{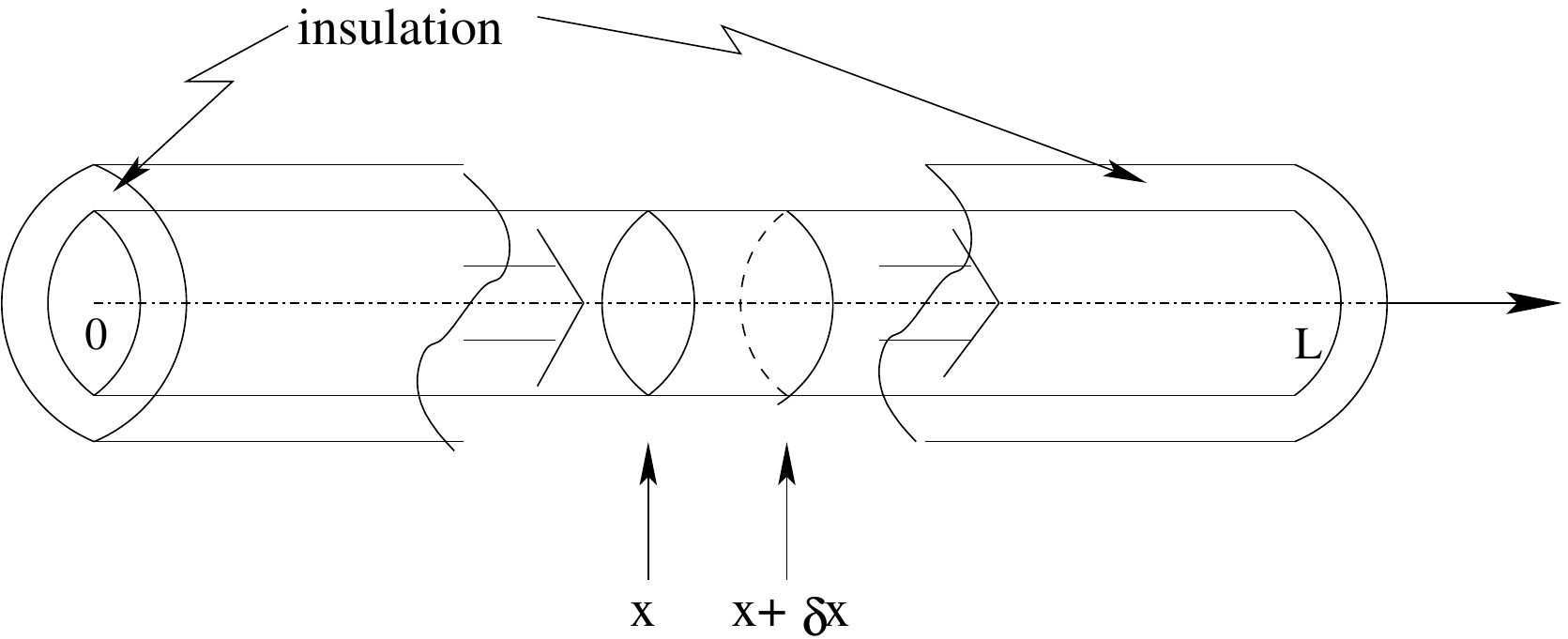}
   \caption{Illustration of a bar with constant cross section, insulated along its lateral surface.}
   \label{barpic}
\end{center}
\end{figure}

We consider a homogeneous bar (constant material parameter throughout the bar), which is insulated along
its lateral surface and has a constant cross-section along its length, see Figure \ref{barpic}. 
Let $u(x , t)$ and $q(x , t)$ denote the temperature and energy flux at cross-section $x$ at time $t$. Focusing
on the segment of the bar between cross-sectional segments $x$ and $x \, + \, \delta x$, applying the
principle of conservation of energy we derive the equation
\be
 \int_{x}^{x \, + \, \delta x} \, c \, \rho \, A \, \frac{\partial u(\xi , t)}{\partial t} \, d\xi
 \ = \ q(x , t) \, A \ - \ q(x \, + \, \delta x) \, A \ + \ 
 \int_{x}^{x \, + \, \delta x} \,  \rho \, A \, \tilde{f}(\xi , t) \, d\xi \, ,
 \label{me1}
\ee
where $c$, $\rho$, and $A$ denote the specific heat constant, the density, and the  
cross-sectional area of the bar, respectively. The function $\tilde{f}$ represents 
an internal energy source. Rearranging \eqref{me1}, using 
\[
q(x , t)  \ - \ q(x \, + \, \delta x) \ = \ - \int_{x}^{x \, + \, \delta x} \frac{\partial q(\xi , t)}{\partial \xi} \, d\xi \, ,
\] 
yields
\be
 \int_{x}^{x \, + \, \delta x} \, \left( c \, \rho \, A \, \frac{\partial u(\xi , t)}{\partial t} \, d\xi
 \ + \  A \, \frac{\partial q(\xi , t)}{\partial \xi} \ - \ 
   \rho \, A \, \tilde{f}(\xi , t)  \right) d\xi  \ = \  0 \, .
 \label{me2}
\ee
With $f(x , t) \, := \, \rho \, \tilde{f}(x , t) $, using the arbitrariness of $x$ and $\delta x$ leads to the pointwise equation
\be
  \frac{\partial u(x , t)}{\partial t}  \ + \ \frac{1}{c \, \rho} \frac{\partial q(x , t)}{\partial x} 
  \ = \  \frac{1}{c \, \rho} f(x , t) \, , \ 0 < x < L \, , \ t > 0 \, .
\label{me2p5}
\ee
The corresponding steady state equation is then
\be
  \frac{d}{dx} q(x) \ = \ f(x) \, ,   \ 0 < x < L \, .
 \label{me3}
\ee
In comparison of \eqref{me3} with  \eqref{reg1} and \eqref{reg2} it is important to note that the ``outside''
derivative comes from the conservation of energy principle and not from the diffusion process.

Fourier's law of thermal conduction (analogous to Fick's law of diffusion) postulates that
\[
 q(x , t) \ = \ - k \, \frac{\partial u(x , t)}{\partial x} \, , \ \ \mbox{ for } k \mbox{ a constant.}
\]
In \cite{sch011} Schumer et al. presented an argument for the diffusion process modeled as a 
random walk process, assuming a ``heavy tail'' distribution for the jumps, to obtain a fractional
law of diffusion given by
\be
  q(x , t)  \ = \ - \kappa \, \left( r \, \mbox{}^{RL}_{0}D_{x}^{(\alpha - 1)}  \ 
  + \ (1 - r) \, \mbox{}^{RL}_{x}D_{L}^{(\alpha - 1)} \right)u(x , t) \, ,
\label{me4}
\ee
for $\kappa$ a dispersion coefficient, $0 < r < 1$ a parameter, and $1 < \alpha < 2$.
With \eqref{me4}, for $r = 1$ the corresponding steady-state diffusion equation with homogeneous 
boundary conditions becomes
\begin{align}
   - \mbox{}^{RL}_{0}D_{x}^{\alpha} u(x) \, := \, - D^{2} \, \bfD^{-(2 - \alpha)} u(x) &= \ f(x) \, , \ \ 0 < x < L \, ,   \label{me5} \\
     u(0) \ = \ u(L) &= 0 \, .   \label{me6}
\end{align}    
We have that the kernel of the operator  $\mbox{}^{RL}_{0}D_{x}^{\alpha}$ is 
$ker( \mbox{}^{RL}_{0}D_{x}^{\alpha} ) \ = \ 
span\{ x^{\alpha - 2} , \, x^{\alpha - 1} \}$.

In sophomore calculus the procedure taught for determining the solution of $2^{nd}$ order, linear
differential equations is to find the general solution to the homogeneous problem (i.e., determine the
kernel of the operator), and then add to it a particular solution. Similarly, we can write the solution
to \eqref{me5} as
\[
  u(x) \ = \ C_{1} x^{\alpha - 2} \ + \ C_{2} x^{\alpha - 1} \ + \ F(x) \, ,
\]
where $F(x)$ satisfies $\mbox{}^{RL}_{0}D_{x}^{\alpha} F(x)  \ = \ - f(x)$. However, because of the singular function
$x^{\alpha - 2}$, the only way for $u$ to satisfy the boundary conditions \eqref{me6} is if $F(x)$ is chosen
to satisfy $F(0) = 0$, e.g. $F(x) \ = \ - \bfD^{- \alpha} f(x)$, and constants $C_{1}$ and $C_{2}$ chosen as
$C_{1} = 0$, $C_{2} = - F(L) /  L^{\alpha - 1}$.

Suppose instead that the boundary conditions associated with \eqref{me5} are $u(0) = u(L) = 1$.
(Thinking of \eqref{me5} as modeling a physical experiment, the change in boundary condition simply amounts
to the engineer performing the experiment relabeling their thermometer by adding a 1 to its values.)
Mathematically, to reduce the new problem to one with homogeneous boundary conditions, we introduce the
change of unknown $v(x) \ = \ u(x) \, - \, 1$, satisfying
\be
     - \mbox{}^{RL}_{0}D_{x}^{\alpha} v(x) \ =  \ - D^{2} \bfD^{-(2 - \alpha)} (u(x) \, - \, 1) 
     \ = \ f(x) \, + \, \frac{1}{\Gamma(1 - \alpha)} x^{-\alpha} \, , \ \ 0 < x < L \, .   \label{me8}
\ee
Now to simulate the problem determined by $v$ requires an infinite energy source be applied at $x = 0$!

In place of  \eqref{me5}, consider the one-sided fractional diffusion equation given by using
\[
   q(x , t) \ = \ - \left(\bfD^{-(2 - \alpha)} \, D \right)u(x , t) \, .
\]
This leads to the steady state equation
\be
  \mcL_{1}^{\alpha} u(x) \, := \, - D \, \bfD^{-(2 - \alpha)} \, D u(x) \ = \ f(x) \, , \ \ 0 < x < L \, .
   \label{me9} 
\ee
The kernel of the operator $\mcL_{1}^{\alpha}$ is $ker(\mcL_{1}^{\alpha}) \, = \, span\{1 , \, x^{\alpha - 1} \}$.
With \eqref{me9} subject to homogeneous boundary conditions, its solution corresponds to that given above.

Again, considering the case of boundary conditions $u(0) = u(L) = 1$. 
Under the change of unknown $v \ = \ u  \, - \, 1$,  \eqref{me9} transforms to
 \be
  \mcL_{1}^{\alpha} v(x) \, := \, - D \, \bfD^{-(2 - \alpha)} \, D v(x) \ = \ f(x) \, , \ \ 0 < x < L \, .
   \label{me10} 
\ee
The simulation of this model equation would require the same energy source as for the case 
$u(0) = u(L) = 0$, which physically makes sense!
For this reason we believe the appropriate interpretation of $\bfD^{\alpha}$ and $\bfD^{\alpha *}$,
$1 < \alpha < 2$, in diffusion
problems is
\be
   \bfD^{\alpha} \, := \, D \, \bfD^{-(2 - \alpha)} \, D \, , \ \  \mbox{ and } \ \
       \bfD^{\alpha *} \, := \, D \, \bfD^{-(2 - \alpha)*} \, D \, .
  \label{me11}
\ee
Note that the definition of  $\bfD^{\alpha}$ and $\bfD^{\alpha *}$
given in \eqref{me11} differs from both the Riemann-Liouville 
and Caputo definitions of $\bfD^{\alpha}$.  \\

\underline{Physical interpretation of the nonlocal diffusion equation} \\
Referring to the setting introduced at the beginning of this section, for a bar at constant
temperature at the atomic scale the particles are in constant motion. However, in relation to the
model depicted in Figure \ref{barpic}, at the macroscopic scale the ``average energy flux'' across 
any cross section is zero.

Consistent with Fourier's law of heat conduction (also Fick's law of diffusion), we posit that a temperature 
gradient across a cross section at location $s$ results in a macroscopic scale having a zero ``average energy flux''
at $s$. If we assume that there is a nonlocal effect from a flux originating at a cross section $s$, proportional
to $1 / \mbox{(distance from that point)}^{(\alpha - 1)}$ then the contribution to the flux at cross section $x$
from points to its left is given by
\be
      k \, \int_{0}^{x} (x \, - \, s)^{(1 - \alpha)} \, (-) \frac{\partial u(s , t)}{\partial s} \, ds \, ,
    \label{me21}
\ee
where $k$ again denotes a thermal conductivity factor.
(The $(-) \frac{\partial u(s , t)}{\partial s}$ denotes the fact that energy flows from ``hot to cold.'')

Similarly, the  contribution to the flux at cross section $x$
from points to its right is given by
\be
     k \,  \int_{x}^{1} (s \, - \, x)^{(1 - \alpha)} \, \frac{\partial u(s , t)}{\partial s} \, ds \, .
    \label{me22}
\ee

Proceeding in an analogous manner to the derivation given at the beginning of this section we obtain 
(corresponding to \eqref{me2p5}) the fractional diffusion equation
\begin{align}
  \frac{\partial u(x , t)}{\partial t}  \ + \ \frac{k}{c \, \rho} 
  \left( \frac{\partial}{\partial x}  \int_{0}^{x} (x \, - \, s)^{(1 - \alpha)} \, (-) \frac{\partial u(s , t)}{\partial s} \, ds \
 \right.  &+ \ \left. \frac{\partial}{\partial x} \int_{x}^{1} (s \, - \, x)^{(1 - \alpha)} \, \frac{\partial u(s , t)}{\partial s} \, ds \right)
  \ = \  f(x , t) \, , \nonumber \\
  & \quad \quad \quad \quad \quad \mbox{ for } \ 0 < x < L \, , \ t > 0 \, ,
\label{me23}
\end{align}
or equivalently written,
\[   \frac{\partial u(x , t)}{\partial t}  \ + \ \frac{k}{c \, \rho} 
  \left( \bfD^{\alpha} u(x , t) \ + \ \bfD^{\alpha *} u(x , t) \right) \ = \ f(x , t) \, , \ 
  \mbox{ for } \ 0 < x < L \, , \ t > 0 .
\]  

\setcounter{equation}{0}
\setcounter{figure}{0}
\setcounter{table}{0}
\setcounter{theorem}{0}
\setcounter{lemma}{0}
\setcounter{corollary}{0}
\section{Kernel of the operator $r \bfD^{\alpha} \, + \, (1 - r) \bfD^{\alpha *}$}
 \label{sec_reg}
In this section we establish the kernel for the operator
\be
 \mcL_{r}^{\alpha} \, = \, - \left( r \bfD^{\alpha} \, + \, (1 - r) \bfD^{\alpha *} \right) \, ,
 \label{eww1}
\ee 
where $ 1 < \alpha < 2$. 

Important in the discussion is the precise definition of the operator (\ref{eww1}). For our interest,
arising from fractional advection-diffusion equations, the operator (\ref{eww1}) is interpreted as
\be
\mcL_{r}^{\alpha} u \ = \ - \left( r \bfD^{\alpha} \, + \, (1 - r) \bfD^{\alpha *} \right) u \ := \ 
-\left( r D \bfD^{-(2 - \alpha)} D  \, + \, (1 - r) D \bfD^{-(2 - \alpha)*} D \right) u \, .
\label{eww2}
\ee

\textbf{Remark}: The definition given in (\ref{eww2}) differs from the Riemann-Liouville definition 
for $\bfD^{\alpha}$, where both integer order derivatives occur after the fractional integral. These different
interpretations represent different operators and hence they have different kernels. For example,
$u = constant$ is in the kernel of the operator defined in (\ref{eww2}). However, $u = constant$
is not in the kernel of (\ref{eww1}) using the Riemann-Liouville definition of the fractional differential
operators.

\subsection{Kernel of  $\mcL_{1/2}^{\alpha}$}
Before discussing the general case we consider the kernel of $D \bfD^{-(2 - \alpha)} D \, + \,
D \bfD^{-(2 - \alpha)*} D$.

\begin{lemma} \label{lma1}
A kernel function of the operator $D \bfD^{-(2 - \alpha)} \, + \,
D \bfD^{-(2 - \alpha)*}$ is
\be
k_{1/2}(x) \ := \, x^{\alpha/2 - 1} (1 - x)^{\alpha/2 - 1} \, .
\label{dk1}
\ee
\end{lemma}
\textbf{Proof}: Using the definition of the fractional integral, we have
\bea
\bfD^{-(2 - \alpha)}k_{1/2}(x) 
&=& \frac{1}{\Gamma(2 - \alpha)} \int_{0}^{x} (x - s)^{1 - \alpha} \, 
              s^{\alpha/2 - 1} \, (1 - s)^{\alpha/2 - 1} \, ds \,  \nonumber \\
 & & \ \ \mbox{(with the substitution $z \, = \, s / x$)}  \nonumber \\
&=& \frac{1}{\Gamma(2 - \alpha)} \, x^{1 - \alpha/2} \int_{0}^{1} (1 - z)^{1 - \alpha} \, z^{\alpha/2 - 1} \,
        (1 \, - \, z x)^{\alpha/2 - 1} \, dz   \nonumber \\
&=& \frac{1}{\Gamma(2 - \alpha)} \, x^{1 - \alpha/2}  
 \frac{\Gamma(\alpha/2) \, \Gamma(2 - \alpha)}{\Gamma(2 - \alpha/2)}  \, 
    \mbox{}_{2}F_{1}(1 - \alpha/2, \, \alpha/2; \, 2 - \alpha/2; \, x)    \nonumber \\
&=&  \frac{\Gamma(\alpha/2)}{\Gamma(2 - \alpha/2)}  \, x^{1 - \alpha/2}   \, 
    \mbox{}_{2}F_{1}(\alpha/2, \, 1 - \alpha/2; \, 2 - \alpha/2; \, x)    \nonumber \\
& & \ \ \mbox{(using  Proposition \ref{HGFp1})}  \nonumber \\
&=& \frac{\Gamma(\alpha/2)}{\Gamma(2 - \alpha/2)}  \, x^{1 - \alpha/2}   \,  
 \frac{\Gamma(2 - \alpha/2)}{\Gamma(1 - \alpha/2) \, \Gamma(1)} 
   \int_{0}^{1} (1 - z)^{0} \, z^{-\alpha/2} \,
        (1 \, - \, z x)^{-\alpha/2} \, dz   \nonumber \\     
&=& \frac{\Gamma(\alpha/2)}{\Gamma(1 - \alpha/2)} \, 
  x^{1 - \alpha/2}  \, x^{-(1 - \alpha/2)} \int_{0}^{x} s^{-\alpha/2} \, (1 - s)^{-\alpha/2} \, ds \, \nonumber \\
&=& \frac{\Gamma(\alpha/2)}{\Gamma(1 - \alpha/2)} \, 
  \int_{0}^{x} s^{-\alpha/2} \, (1 - s)^{-\alpha/2} \, ds \, .
  \label{eww5}
\eea  
Then, from (\ref{eww5}),
\be
D \bfD^{-(2 - \alpha)} k_{1/2}(x) \ = \ \frac{\Gamma(\alpha/2)}{\Gamma(1 - \alpha/2)}  \, 
  x^{-\alpha/2} \, (1 - x)^{-\alpha/2} \, .
 \label{eww6}
\ee

Next, 
\bea
\bfD^{-(2 - \alpha)*}k_{1/2}(x) 
&=& \frac{1}{\Gamma(2 - \alpha)} \int_{x}^{1} (s - x)^{1 - \alpha} \, 
              s^{\alpha/2 - 1} \, (1 - s)^{\alpha/2 - 1} \, ds \,  \nonumber \\
 & & \ \ \mbox{(with the substitution $z \, = \, (1 - s) /(1 - x)$)}  \nonumber \\
&=& \frac{1}{\Gamma(2 - \alpha)} \, (1 - x)^{1 - \alpha/2} \int_{0}^{1} (1 - z)^{1 - \alpha} \, z^{\alpha/2 - 1} \,
        (1 \, - \, z (1 - x) )^{\alpha/2 - 1} \, dz   \nonumber \\
&=& \frac{1}{\Gamma(2 - \alpha)} \, (1 - x)^{1 - \alpha/2}  \, 
 \frac{\Gamma(\alpha/2) \, \Gamma(2 - \alpha)}{\Gamma(2 - \alpha/2)}  \, 
    \mbox{}_{2}F_{1}(1 - \alpha/2, \, \alpha/2; \, 2 - \alpha/2; \, 1 - x )    \nonumber \\
&=&  \frac{\Gamma(\alpha/2)}{\Gamma(2 - \alpha/2)} \, (1 - x)^{1 - \alpha/2}    \, 
    \mbox{}_{2}F_{1}(\alpha/2, \, 1 - \alpha/2; \, 2 - \alpha/2; \, 1 - x )    \nonumber \\
& & \ \ \mbox{(using  Proposition \ref{HGFp1})}  \nonumber \\
&=&\frac{\Gamma(\alpha/2)}{\Gamma(2 - \alpha/2)} \, (1 - x)^{1 - \alpha/2}    \, 
 \frac{\Gamma(2 - \alpha/2)} {\Gamma(1 - \alpha/2) \, \Gamma(1)}
   \int_{0}^{1} (1 - z)^{0} \, z^{-\alpha/2} \,
        (1 \, - \, z (1 - x) )^{-\alpha/2} \, dz   \nonumber \\     
&=& \frac{\Gamma(\alpha/2)}{\Gamma(1 - \alpha/2)} \, 
  (1 - x)^{1 - \alpha/2}  \, (1 - x)^{-(1 - \alpha/2)} \int_{x}^{1} s^{-\alpha/2} \, (1 - s)^{-\alpha/2} \, ds \, \nonumber \\
&=& \frac{\Gamma(\alpha/2)}{\Gamma(1 - \alpha/2)} \,
  \int_{x}^{1} s^{-\alpha/2} \, (1 - s)^{-\alpha/2} \, ds \, .
  \label{eww7}
\eea  
Then, from (\ref{eww7}),
\be
D \bfD^{-(2 - \alpha)*} k_{1/2}(x) \ = \ - \frac{\Gamma(\alpha/2)}{\Gamma(1 - \alpha/2)} \, 
  x^{-\alpha/2} \, (1 - x)^{-\alpha/2} \, .
 \label{eww8}
\ee
Combining (\ref{eww6}) and (\ref{eww8}) we obtain
$(D \bfD^{-(2 - \alpha)} \, + \,D \bfD^{-(2 - \alpha)*}) k_{1/2}(x) \ = \ 0$.  \\
\mbox{ } \hfill \qed

\begin{figure}[!ht]  
\begin{center}
 \includegraphics[height=2.5in]{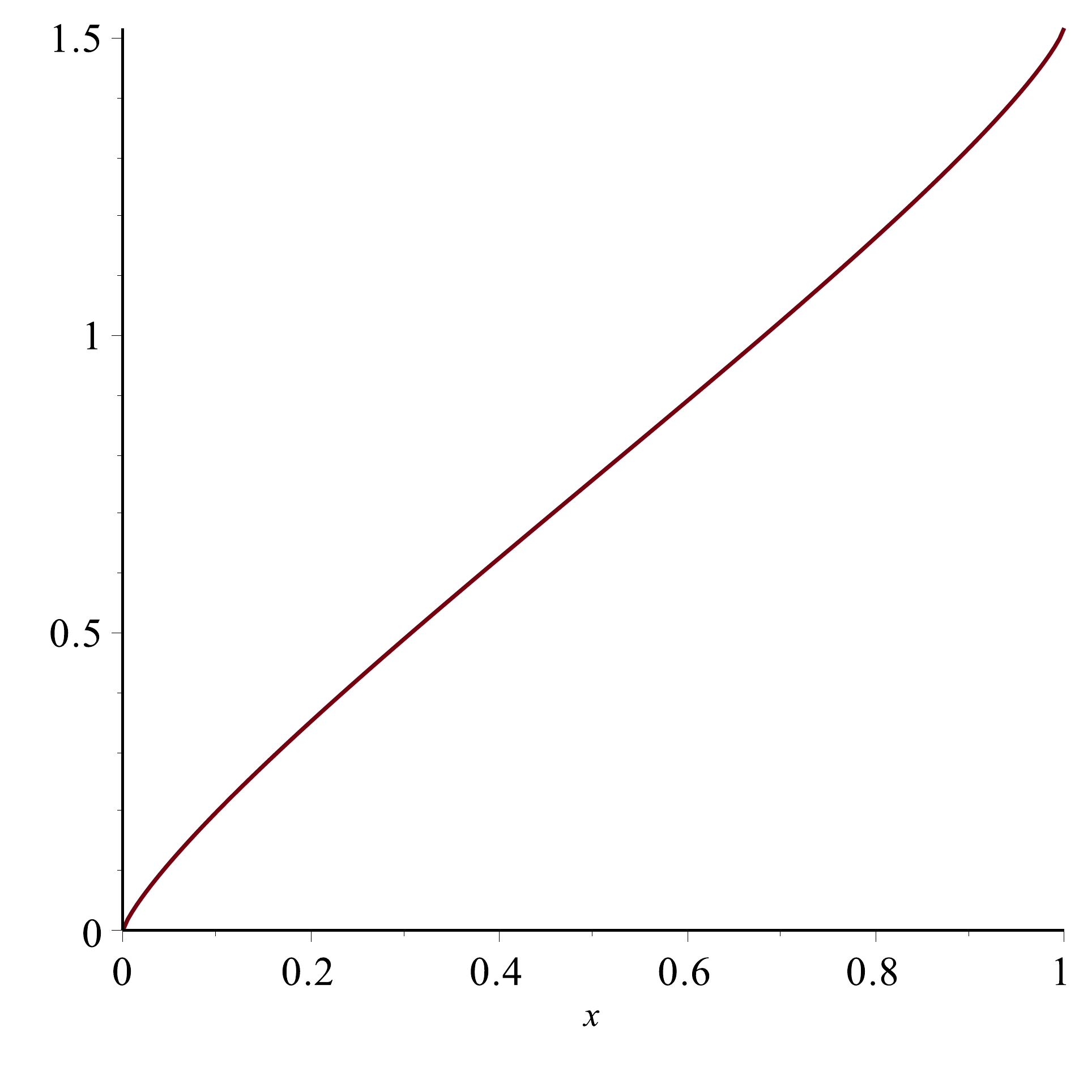}
   \caption{Plot of $K_{1/2}(x)$ for $\alpha = 1.6$.}
   \label{figK1ov2}
\end{center}
\end{figure}

Let 
\be
K_{1/2}(x) := \int_{0}^{x} k_{1/2}(s) \, ds \ = \ \frac{2}{\alpha} x^{\alpha/2} \mbox{}_{2}F_{1}(\alpha/2 , 1 - \alpha/2 ; 
\, 1 + \alpha/2 ; \, x) \, .
\label{defK1}
\ee
A plot of $K_{1/2}(x)$ for $\alpha = 1.6$ is given in Figure \ref{figK1ov2}.

\begin{lemma} \label{ker1}
The kernel of $\mcL_{1/2}^{\alpha}(\cdot)$, $ker(\mcL_{1/2}^{\alpha})$, is given by
$ker(\mcL_{1/2}^{\alpha}) \ = \ span\{1 , K_{1/2}(x)\}$.
\end{lemma}
\textbf{Proof}: From above it is clear that $span\{1 , K_{1/2}(x)\} \subset ker(\mcL_{1/2}^{\alpha})$. What remains is to
show that $dim(ker(\mcL_{1/2}^{\alpha})) \, = \, 2$. 

With $z(x) \, = \, 1 + x$ and $f(x) \ = \ \frac{-1}{2} \, \frac{1}{\Gamma(2 - \alpha)} \, x^{1 - \alpha} 
            \ + \  \frac{1}{2} \, \frac{1}{\Gamma(2 - \alpha)} \, (1 - x)^{1 - \alpha}$, 
a straightforward calculation shows that  $\mcL_{1/2}^{\alpha} z(x) \, = \, f(x)$ on $I$.      
As $K_{1/2}(1) \ne 0$ we can choose $c_{1}$ and $c_{2}$
such that $\hat{z}(x) \ := \ z(x) \, + \, c_{1} 1 \, + \, c_{2} K_{1/2}(x)$ satisfies $\hat{z}(0) = \hat{z}(1) = 0$ and
$\mcL_{1/2}^{\alpha} \hat{z}(x) \, = \, f(x)$. 

Suppose there was another linearly independent function $s(x) \in ker(\mcL_{1/2}^{\alpha})$. Without loss of generality
we may assume that $s(0) = s(1) = 0$. (If this was not the case we would form a linear combination of $s(x)$ with the
other two linearly independent kernel function $1$ and $K_{1/2}(x)$.) Then $\tilde{z}(x) \ := \ \hat{z}(x) \, + \, s(x)$ satisfies
$\tilde{z}(0) = \tilde{z}(1) = 0$ and $\mcL_{1/2}^{\alpha} \tilde{z}(x) \, = \, f(x)$. However, the existence of
$\tilde{z}(x) \ne \hat{z}(x)$ contradicts the uniqueness of the solution 
 to $\mcL_{1/2}^{\alpha} u(x)  \, = \, f(x)$, with  $u(0) = u(1) = 0$, \cite{erv061}. 
\\
\mbox{ } \hfill \qed

\subsection{Kernel of  $\mcL_{r}^{\alpha}$}
In the section we extend the discussion from the previous section to the operator
\be
\mcL_{r}^{\alpha} u \ = \ 
- \left( r D \bfD^{-(2 - \alpha)} D  \, + \, (1 - r) D \bfD^{-(2 - \alpha)*} D \right) u \, .
\label{ewq1}
\ee

\begin{lemma} \label{ker1r}
With 
$k(x) \ := \ x^{p} \, (1 - x)^{q}$, $K(x) \ := \ \int_{0}^{x} k(s) \, ds$, we have that $K(x) \in ker(\mcL_{r}^{\alpha})$ if
\begin{align}
  (i)   \quad 3 - \alpha + p + q &= 1 \, ,  \label{propK1} \\
 \mbox{and } \   (ii)  \quad  r \sin( \pi (-q) ) &= (1 - r) \sin( \pi (-p) ) \, .  \label{propK2}
\end{align}
\end{lemma}
\textbf{Proof}: Proceeding as above,
\bea
\bfD^{-(2 - \alpha)}k(x) 
&=& \frac{1}{\Gamma(2 - \alpha)} \int_{0}^{x} (x - s)^{1 - \alpha} \, 
              s^{p} \, (1 - s)^{q} \, ds \,  \nonumber \\
&=& \frac{1}{\Gamma(2 - \alpha)} \, x^{2 - \alpha + p} \int_{0}^{1} (1 - z)^{1 - \alpha} \, z^{p} \,
        (1 \, - \, z x)^{q} \, dz   \ \ \mbox{ (using $z = s / x$)} \nonumber \\
&=& \frac{1}{\Gamma(2 - \alpha)} \, x^{2 - \alpha + p}  
 \frac{\Gamma(p + 1) \, \Gamma(2 - \alpha)}{\Gamma(3 - \alpha + p)}  \, 
    \mbox{}_{2}F_{1}(-q, \, p + 1; \, 3 - \alpha + p; \, x)    \nonumber \\
& & \mbox{(provided $3 - \alpha + p > p + 1 > 0$, which is true for $1 < \alpha < 2$)}  \nonumber \\    
&=&  \frac{\Gamma(p + 1)} {\Gamma(3 - \alpha + p)}\, x^{2 - \alpha + p}    \, 
    \mbox{}_{2}F_{1}(p + 1, \, -q ; \, 3 - \alpha + p; \, x)    \nonumber \\
& & \ \ \mbox{(using  Proposition \ref{HGFp1}, provided $3 - \alpha + p > -q > 0$)}  \nonumber \\
&=& \frac{\Gamma(p + 1)} {\Gamma(3 - \alpha + p)}\, x^{2 - \alpha + p}    \,   \cdot \nonumber \\
 & & \ \ \ \ \ \ \ \ \ \ 
 \frac{\Gamma(3 - \alpha + p)}{\Gamma(-q) \, \Gamma(3 - \alpha + p + q)}
   \int_{0}^{1} (1 - z)^{2 - \alpha + p + q} \, z^{-q - 1} \,
        (1 \, - \, z x)^{-p - 1} \, dz   \nonumber \\     
&=& \frac{\Gamma(p + 1)}{\Gamma(-q) \, \Gamma(3 - \alpha + p + q)} \, 
  x^{2 - \alpha + p}   \, x^{-(2 - \alpha + p)}  
  \int_{0}^{x} (x - s)^{2 - \alpha + p + q} \, s^{-q - 1} \,
        (1 \, - \, s)^{-p - 1} \, ds \, \nonumber \\
&=& \frac{\Gamma(p + 1)}{\Gamma(-q) \, \Gamma(3 - \alpha + p + q)} \, 
  \int_{0}^{x} (x - s)^{2 - \alpha + p + q} \, s^{-q - 1} \,
        (1 \, - \, s)^{-p - 1} \, ds  \, \nonumber \\
&=&     \frac{\Gamma(p + 1)}{\Gamma(-q)} \, 
        \bfD^{-(3 - \alpha + p + q)} x^{-q - 1} \,  (1 \, - \, x)^{-p - 1} \, .
  \label{ewq2}
\eea  

Next,
\bea
\bfD^{-(2 - \alpha)*}k(x) 
&=& \frac{1}{\Gamma(2 - \alpha)} \int_{x}^{1} (s - x)^{1 - \alpha} \, 
              s^{p} \, (1 - s)^{q} \, ds \,  \nonumber \\
&=& \frac{1}{\Gamma(2 - \alpha)} \, (1 - x)^{2 - \alpha + q} \int_{0}^{1} (1 - z)^{1 - \alpha} \, z^{q} \,
        (1 \, - \, z (1 - x) )^{p} \, dz  \nonumber \\
& &          \ \mbox{ (using $z = (1 - s)/(1 - x)$)}  \nonumber \\
&=& \frac{1}{\Gamma(2 - \alpha)} \, (1 - x)^{2 - \alpha + q}  
 \frac{\Gamma(q + 1) \, \Gamma(2 - \alpha)}{\Gamma(3 - \alpha + q)}  \, 
    \mbox{}_{2}F_{1}(-p, \, q + 1; \, 3 - \alpha + q; \, (1 - x) )    \nonumber \\
& & \mbox{(provided  $3 - \alpha + q > q + 1 > 0$, which is true for $1 < \alpha < 2$)} \nonumber \\
&=&  \frac{\Gamma(q + 1)}{\Gamma(3 - \alpha + q)}  \, (1 - x)^{2 - \alpha + q}   \, 
    \mbox{}_{2}F_{1}(q + 1, \, -p ; \, 3 - \alpha + q; \, (1 - x) )    \nonumber \\
& & \ \ \mbox{(using  Proposition \ref{HGFp1}, provided $3 - \alpha + q > -p > 0$)}  \nonumber \\
&=&    \frac{\Gamma(q + 1)}{\Gamma(3 - \alpha + q)}  \, (1 - x)^{2 - \alpha + q}   \,    \cdot \nonumber \\
 & & \ \ \ \ \ \ \ \ \ \ 
 \frac{\Gamma(3 - \alpha + q)}{\Gamma(-p) \, \Gamma(3 - \alpha + p + q)} 
   \int_{0}^{1} (1 - z)^{2 - \alpha + p + q} \, z^{-p - 1} \,
        (1 \, - \, z (1 - x) )^{-q - 1} \, dz   \nonumber \\     
&=&  \frac{\Gamma(q + 1)}{\Gamma(-p) \, \Gamma(3 - \alpha + p + q)}  \, (1 - x)^{2 - \alpha + q}   \,  
  (1 - x)^{-(2 - \alpha + q)} \cdot \nonumber \\
  & &  \ \ \ \ \ \ \ \ \ \ 
  \int_{x}^{1} (s - x)^{2 - \alpha + p + q} \, s^{-q - 1} \, (1 - s)^{-p - 1} \, ds \, \nonumber \\
&=&  \frac{\Gamma(q + 1)}{\Gamma(-p) \, \Gamma(3 - \alpha + p + q)} \, 
  \int_{x}^{1} (s - x)^{2 - \alpha + p + q} \, s^{-q - 1} \, (1 - s)^{-p - 1} \, ds \, \nonumber \\
&=&   \frac{\Gamma(q + 1)}{\Gamma(-p)}  \, 
  \bfD^{-(3 - \alpha + p + q)*} x^{-q - 1} \, (1 - x)^{-p - 1} \, .
  \label{ewq3}
\eea  

Comparing \eqref{ewq2} and \eqref{ewq3}, 
$r \, D \, D^{-(2 - \alpha)} k(x) \ + \ (1 - r) \, D \, D^{-(2 - \alpha) *} k(x) \ = \ 0$ if 
\bea
(i) \ \ \ 3 - \alpha + p + q &=&1 \, ,   \label{ewq4}  \\
(ii) \ \ \  r \frac{\Gamma(p + 1)}{\Gamma(-q) \, \Gamma(3 - \alpha + p + q)}
 &=& (1 - r) \frac{\Gamma(q + 1)}{\Gamma(-p)\, \Gamma(3 - \alpha + p + q)}  \nonumber \\
 \Longleftrightarrow \ \ 
    r \frac{\Gamma(p + 1)}{\Gamma(-q)}
 &=& (1 - r) \frac{\Gamma(q + 1)}{\Gamma(-p)}  \nonumber  \\
 \Longleftrightarrow \ \ 
 r \, \Gamma(-p) \, \Gamma(1 - (-p)) & = & (1 - r) \,  \Gamma(-q) \, \Gamma(1 - (-q)) \nonumber \\
 \Longleftrightarrow \ \ 
 r \frac{\pi}{\sin(\pi (-p))}  &=& (1 - r)  \frac{\pi}{\sin(\pi (-q))} \nonumber \\
  \mbox{(using  } \Gamma(1 - z) \, \Gamma(z) \ &=& \ \pi / \sin( \pi z) \, , \mbox{ valid for } z \neq 0, \, \pm 1, \, \pm 2, 
  \ldots   )  \nonumber \\
  \Longleftrightarrow \ \ 
 r \sin(\pi \, (-q))  &=& (1 - r) \sin(\pi \, (-p)) \, .  \label{ewq5}
\eea
\mbox{ } \hfill \qed

\begin{corollary} \label{cor1r}
The kernel of $\mcL_{r}^{\alpha}(\cdot)$, $ker(\mcL_{r}^{\alpha})$, is given by
$ker(\mcL_{1/2}^{\alpha}) \ = \ span\{1 , K(x)\}$, where $K(x)$, given in Lemma \ref{ker1r}, may be written 
as $K(x) \ = \ \int_{0}^{x} k(s) \, ds \ = \ \frac{1}{p + 1} x^{p + 1} \, \mbox{}_{2}F_{1}(-q , p + 1 \, ; \, p + 2 \, ; \, x)$. 
\end{corollary}
\mbox{ } \hfill \qed

\underline{Example 2.1}.  The case $r = 1/2$. This corresponds to $\mcL^{\alpha}_{1/2}$. \\
For $r = 1/2$, from (\ref{ewq5}), $p = q$. Then, using (\ref{ewq4}), we have $p \ = \ q \ = \alpha/2 \, - \, 1$, which 
agrees with $k_{1/2}(x)$ given in (\ref{dk1}).

\underline{Example 2.2}.  The case $r \rightarrow 1$. This corresponds to $\mcL^{\alpha}_{1}(u) \ = \
- D \bfD^{-(2 - \alpha)} D(u)$. For this case the kernel is $span\{1 , x^{\alpha - 1}\}$. \\
Now, from (\ref{ewq5}), as $r \rightarrow 1$ then 
\[
  \sin(\pi (-q)) \rightarrow 0 \ \Longrightarrow \ q \rightarrow 0 \, . 
\]
Hence from (\ref{ewq4}) $p \rightarrow \alpha - 2  \
 \Longrightarrow \ K(x) \ = \ x^{\alpha - 1}$ .  \\

\begin{figure}[!ht]
\begin{minipage}{.46\linewidth}

\begin{center}
 \includegraphics[height=2.5in]{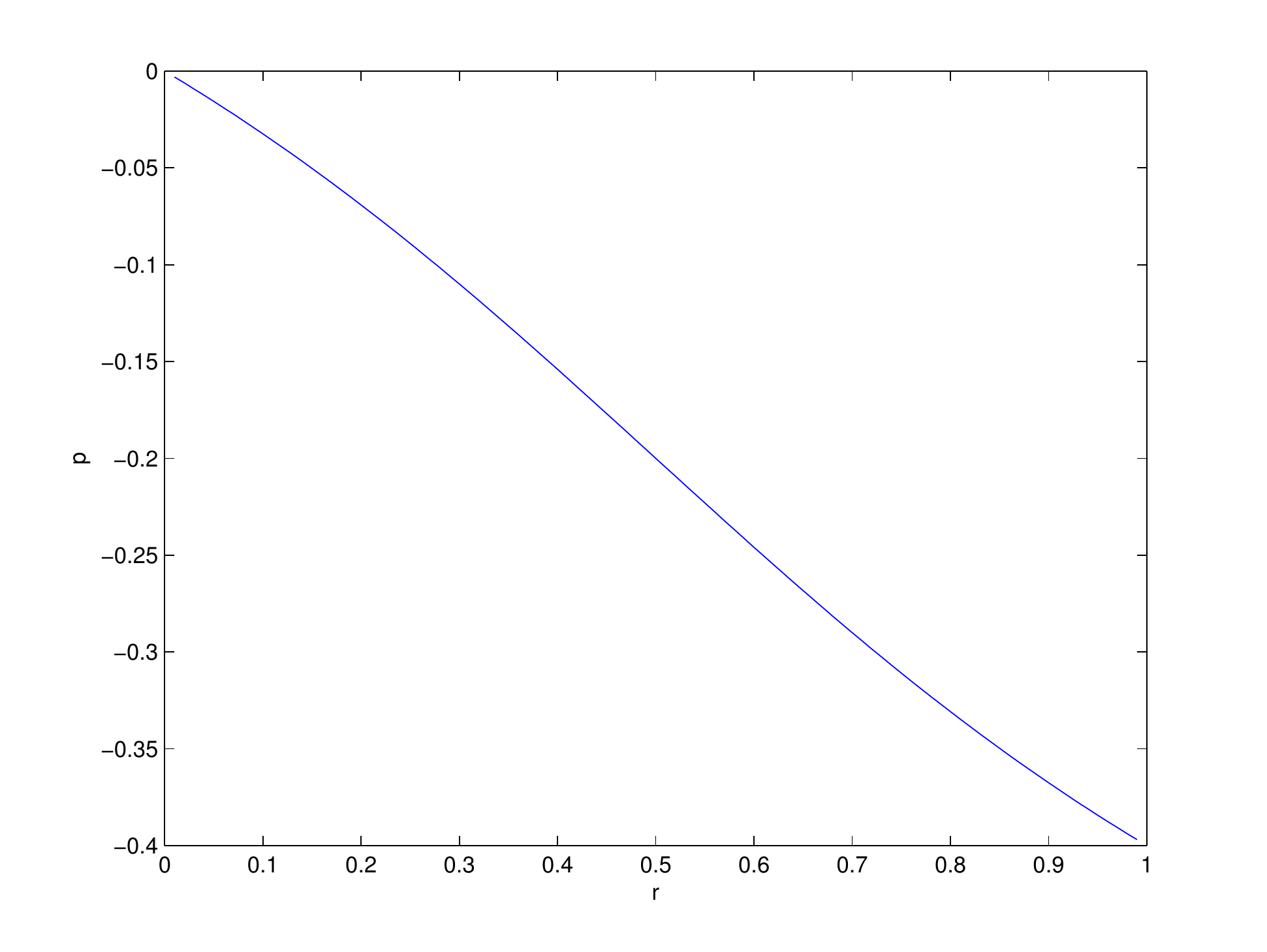}
   \caption{p values solving (\ref{ewq4}) and (\ref{ewq5}) for $\alpha = 1.6$.}
   \label{pvals1}
\end{center}

\end{minipage} \hfill
\begin{minipage}{.46\linewidth}
 
\begin{center}
\includegraphics[height=2.5in]{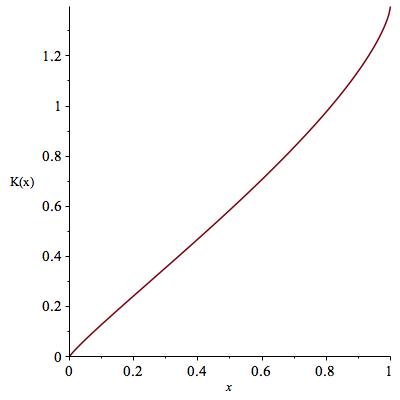}
   \caption{Plot of $K(x)$ for $\alpha = 1.6$ and $r = 0.2764$ (i.e., $p = -0.1$, $q = -0.3$).}
   \label{plotK}
\end{center}
\end{minipage} 
\end{figure}

\begin{lemma} \label{lmaK1}
For $1 \le \alpha < 1.5$, $D \bfD^{-(2 - \alpha)} D$ maps from $H^{\alpha}(I)$ onto $L^{2}(I)$.
\end{lemma}
\textbf{Proof}: We have that $D \, : \, H^{\alpha}(I) \longrightarrow H^{\alpha - 1}(I)$. Now, for $1 < \alpha < 1.5$, 
then $0 < \alpha - 1 < 0.5$, hence $H^{\alpha - 1}(I) \, = \, \tilde{H}^{\alpha - 1}(I)$. As
$D \bfD^{-(2 - \alpha)} \, = \, \mbox{}^{RL}_{0}D_{x}^{\alpha - 1}$, then from Theorem 3.1 \cite{jin151}
$D \bfD^{-(2 - \alpha)} \, : \, \tilde{H}^{\alpha - 1}(I) \longrightarrow L^{2}(I)$. 

To establish that the mapping is onto, we have that for $f \in L^{2}(I)$, 
$D \bfD^{-(2 - \alpha)} D \, u \ = \ f$, where $u \ = \ \frac{1}{\Gamma(\alpha)} 
\int_{0}^{x} (x \, - \, s)^{\alpha - 1} \, f(s) \, ds \, .$ \\
\mbox{ } \hfill \qed

\begin{corollary} \label{corK1}
For $1 \le \alpha < 1.5$, $r \in \real$, $\mcL_{r}^{\alpha}$ maps from $H^{\alpha}(I)$ \underline{into} $L^{2}(I)$.
\end{corollary}
\textbf{Proof}: An analogous argument to that given in the proof of Lemma \ref{lmaK1} establishes that
$D \bfD^{-(2 - \alpha) *} D$ maps from $H^{\alpha}(I)$ onto $L^{2}(I)$. The stated result then follows. \\
\mbox{ } \hfill \qed

In order to give a concise description of the range of $\mcL_{r}^{\alpha}$, with domain $H^{\alpha}(I)$, let
\[
  X^{(1 - \alpha)} \ := \ \{ f \, : \, f(x) \, = \, c x^{1 - \alpha} , \, c \in \real \} \quad \mbox{and} \quad
  X^{(1 - \alpha)*} \ := \ \{ f \, : \, f(x) \, = \, c (1 - x)^{1 - \alpha} , \, c \in \real \} \, .
\]  

\begin{lemma} \label{lmamap2}
For $1 \le \alpha < 2$ $\mcL_{r}^{\alpha}$ maps from $H^{\alpha}(I)$  \underline{into} 
$L^{2}(I) \oplus X^{(1 - \alpha)} \oplus X^{(1 - \alpha)*}$.
\end{lemma}
\textbf{Proof}: The case for $1 \le \alpha < 1.5$ is covered by Corollary \ref{corK1}. For $f(x) \in H^{\alpha}(I)$,
$\alpha \ge 1.5$, let $p(x)$ denote the Hermite cubic interpolating polynomial of $f(x)$. Namely,
\begin{align*}
p(x) &= \ (2 x^{3} \, - \, 3 x^{2} \, + \, 1) f(0) \ + \ (x^{3} \, - \, 2 x^{2} \, + \, x) f'(0) \ + \ 
  (-2 x^{3} \, + \, 3 x^{2}) f(1) \, + \, (x^{3} \, - \, x^{2}) f'(1)   \\
  &= (-2 (1 - x)^{3} \, + \, 3 (1 - x)^{2} ) f(0) \ + \ (-(1 - x)^{3} \, + \, (1 - x)^{2}) f'(0)  \\
  & \quad \ + \
       (2 (1 - x)^{3} \, - \, 3 (1 - x)^{2} \, + \, 1) f(1) \ + \ (-(1 - x)^{3} \, + \, 2 (1 - x)^{2} \, - \, (1 - x)) f'(1) \, .
\end{align*}        
Also, let $\tilde{f}(x) \ = \ f(x) \, - \, p(x) \, \in \tilde{H}^{\alpha}(I)$. From Theorem 2.1 \cite{jin151}, 
$\mcL_{r}^{\alpha} \tilde{f}(x) \in L^{2}(I)$.

Now,
\begin{align*}
\mcL_{r}^{\alpha} f(x) &= \ \mcL_{r}^{\alpha} \tilde{f}(x) \ + \ r \left( f(0) \right) \bfD^{\alpha} 1 \ + \
  r \left( f'(0) \right) \bfD^{\alpha} x \ + \ r \left( -3 f(0) \, - \, 2 f'(0) \, + \, 3 f(1) \, - \, f'(1) \right) \bfD^{\alpha} x^{2} \\
  & \quad 
   \ + \ r \left( 2 f(0) \, + \, f'(0) \, - \, 2 f(1) \ + \ f'(1) \right) \bfD^{\alpha} x^{3}  
\ + \ 
 (1 - r) \left( f(1) \right) \bfD^{\alpha *} 1 \\
& \quad  \ + \  (1 - r) \left( -f'(1) \right) \bfD^{\alpha *} (1 - x) \ 
  + \ (1 - r) \left( -3 f(1) \, + \, 2 f'(1) \, + \, 3 f(0) \, + \, f'(0) \right) \bfD^{\alpha *} (1 - x)^{2} \\
& \quad  \ + \  (1 - r) \left( 2 f(1) \, - \, f'(1) \, - \, 2 f(0) \ - \ f'(0) \right) \bfD^{\alpha *} (1 - x)^{3} \, .
\end{align*}
As $\bfD^{\alpha} 1 \ = \ \bfD^{\alpha *} 1 \ = \ 0$;   $\bfD^{\alpha} x^{2}$, $\bfD^{\alpha} x^{3}$,
$\bfD^{\alpha *} (1 - x)^{2}$, $\bfD^{\alpha *} (1 - x)^{3}  \in L^{2}(I)$, the stated result follows. \\
\mbox{} \hfill \qed

\setcounter{equation}{0}
\setcounter{figure}{0}
\setcounter{table}{0}
\setcounter{theorem}{0}
\setcounter{lemma}{0}
\setcounter{corollary}{0}
\section{Convergence of the Finite Element Method Approximation} 
 \label{secFEMcvg}
In a finite element method (FEM) approximation to \eqref{reg2}
the regularity of the solution $u$ plays a fundamental role in the rate of convergence of the approximation
$u_{h}$ to $u$. In this section we present four numerical experiments and compare the numerical rate
of convergence of the FEM approximation to that predicted theoretically.

From \cite{erv061}, with $X \ = \ \tilde{H}^{\alpha/2}(I)$, the weak formulation of \eqref{reg2} is: 
\textit{Given $f \in H^{-\alpha/2}(I)$ determine $u \in X$ satisfying}
\be
B(u \, , \,  v) \ = \ \langle f \, , \, v \rangle \, , \forall v \in X \, ,
\label{treg1}
\ee
where, $\langle \cdot \, , \, \cdot \rangle$ denotes the $L^{2}$ duality pairing between $H^{-\alpha/2}(I)$
and $\tilde{H}^{\alpha/2}(I)$, and $B(\cdot , \cdot) \, : \, X \times X \longrightarrow \real$ is defined by
\be
 B(w , v) \ := \ r \left(D^{-(2 - \alpha)/2} D w \, , \, D^{-(2 - \alpha)/2 *} D v \right) \ + \
 (1 - r) \left(D^{-(2 - \alpha)/2 *} D w \, , \, D^{-(2 - \alpha)/2} D v \right) \, .
 \label{treg1p5}
\ee 

For $0 = x_{0} < x_{1} < \ldots < x_{N} = 1$ denoting a quasi-uniform partition of $I := (0 , 1)$, $X_{h} \subset X$
denoting the space of continuous, piecewise polynomials of degree $\le k$ on the partition, the finite
element approximation $u_{h} \in X_{h}$ to $u$ is given by
\be
B(u_{h} \, , \,  v_{h}) \ = \ \langle f \, , \, v_{h} \rangle \, , \forall v_{h} \in X_{h} \, .
\label{treg2}
\ee
  
Assuming that $f$ is sufficiently regular such that the regularity of $u$ is determined by the 
kernel of $\mcL_{r}^{\alpha}$, we have the following a priori error bounds, for $C > 0$ a constant
and any $\epsilon > 0$ and $\delta > 0$.
\begin{align}
\| u \, - \, u_{h} \|_{\tilde{H}^{\alpha/2}} &\le  C \, \inf_{v_{h} \in X_{h}}
\| u \, - \, v_{h} \|_{\tilde{H}^{\alpha/2}} \nonumber \\
&\le C \, \left\{ \begin{array}{rl}
  h^{1/2 - \eps} \| u \|_{H^{\alpha/2 + 1/2 - \eps}} \, , & r = 1/2 \, ,   \\
  h^{\min\{p , q\}  + 3/2 - \alpha/2 - \eps} \| u \|_{H^{\min\{p , q\}  + 3/2  - \eps}} \, , & r \ne 1/2 \, ,
  \end{array} \right.     \label{treg3}
\end{align}
where $p$ and $q$ satisfy \eqref{ewq4} and \eqref{ewq5}. 

An application of the Aubin-Nitsche trick yields the following $L^{2}$ a priori error bounds.
\be
\| u \, - \, u_{h} \| \ \le  \ C \, \left\{ \begin{array}{rl}
  h^{1 \, - \, 2 \eps} \| u \|_{H^{\alpha/2 + 1/2 - \eps}} \, , & r = 1/2 \, ,   \\
  h^{2(\min\{p , q\}  + 3/2 - \alpha/2) \, - \, 2 \eps} \| u \|_{H^{\min\{p , q\}  + 3/2  - \eps}} \, , & r \ne 1/2 \, .
  \end{array} \right.     \label{treg4}
\ee
For the Aubin-Nitsche trick the regularity of the associated adjoint problem is the same as that for
$u$ (assuming $f \in L^{2}(I)$). Hence the $L^{2}$ a priori error bound is simply twice that for $H^{\alpha/2}$.

For Examples 1 and 2 the true solution $u$ was chosen to be $x \, + \, kerfun(x)$, with 
$kerfun(x) \in  ker(\mcL_{r}^{\alpha})$ chosen such that $u$ satisfies $u(0) = u(1) = 0$. In Examples 3 and
4 the right hand side $f(x)$ was chosen to be a constant. Results are reported for $\alpha = 1.4$ and $\alpha = 1.6$.
Computations were also performed for $\alpha = 1.2$ and $\alpha = 1.8$ (not included) which exhibited similar
behavior. The $| u - u_{h}|_{H^{\alpha/2}}$ data presented in the tables denotes the 
Slobodetskii semi-norm. The approximation space $X_{h}$ used was the continuous, affine functions on a
uniform partition of $I$.

\underline{Example 1}. \\
With $\alpha = 1.4$, $r = 1/2$, 
\be
 u(x) \ = \ x \ - \ C \,  x^{\alpha/2} \, \mbox{}_{2}F_{1}(\alpha/2 \, , \, 1 - \alpha/2  \, ; \, 1 + \alpha/2 \, , \, x) \, ,
   \label{ppo10}
\ee
where $C \, = \,  \left(  \mbox{}_{2}F_{1}(\alpha/2 \, , \, 1 - \alpha/2  \, ; \, 1 + \alpha/2 \, , \, 1) \right)^{-1} $. \\
The corresponding right hand side is
\be
  f(x) 
     \ =  \ \frac{-1}{2} \, \frac{1}{\Gamma(2 - \alpha)} \, x^{1 - \alpha} 
            \ + \  \frac{1}{2} \, \frac{1}{\Gamma(2 - \alpha)} \, (1 - x)^{1 - \alpha} \, .
            \label{ppo11}
\ee    
The numerical results are presented in Table \ref{Table1}.

\begin{table}[htb!]
\centering
\begin{tabular}{|c||c|c||c|c|} \hline
     $h$ & $ | u  -  u_{h} |_{H^{\alpha/2}(I)}$ & Cvg. rate &
	   $\| u  -  u_{h} \|_{L^{2}(I)}$ & Cvg. rate  \\
\hline
1/64  &  4.209E-02  &    &  8.402E-04  &       \\
1/128  &  2.962E-02  &  0.51  &  4.016E-04  &  1.07     \\
1/256  &  2.088E-02  &  0.50  &  1.936E-04  &  1.05    \\
1/512  &  1.475E-02  &  0.50  &  9.407E-05  &  1.04   \\
1/1024  &  1.042E-02  &  0.50  &  4.598E-05  &  1.03   \\
1/2048  &  7.364E-03  &  0.50  &  2.258E-05  &  1.03  \\  \hline
Pred. &  & 0.50 &   & 1.0 \\
\hline
\end{tabular}
\caption{Example 1. Convergence rates for $\alpha = 1.4$ and $r = 1/2$.}
\label{Table1}
\end{table}

\underline{Example 2}. \\
With $\alpha = 1.4$, $p = -0.15$, $q = \alpha - p - 2$, $r  \, = \,  \sin(\pi p) / ( \sin(\pi p) + \sin(\pi q) )$
\be
 u(x) \ = \ x \ - \ C \,  x^{(p + 1)} \, \mbox{}_{2}F_{1}(-q \, , \, p + 1  \, ; \, p + 2 \, , \, x) \, ,
   \label{ppo12}
\ee
where $C \, = \,  \left(  \mbox{}_{2}F_{1}(-q \, , \, p + 1  \, ; \, p + 2 \, , \, 1) \right)^{-1} $. \\
The corresponding right hand side is
\begin{align}
  f(x) 
     &=   -r \, \frac{1}{\Gamma(2 - \alpha)} \, x^{1 - \alpha} 
            \ + \  (1 - r) \, \frac{1}{\Gamma(2 - \alpha)} \, (1 - x)^{1 - \alpha} \, .
            \label{ppo13}
\end{align}            
The numerical results are presented in Table \ref{Table2}.

\begin{table}[htb!]
\centering
\begin{tabular}{|c||c|c||c|c|} \hline
     $h$ & $ | u  -  u_{h} |_{H^{\alpha/2}(I)}$ & Cvg. rate &
	   $\| u  -  u_{h} \|_{L^{2}(I)}$ & Cvg. rate  \\
\hline
1/64  &  1.463E-01  &    &  1.609E-03  &     \\  
1/128  &  1.146E-01  &  0.35  &  7.847E-04  &  1.04    \\ 
1/256  &  8.990E-02  &  0.35  &  3.831E-04  &  1.03    \\ 
1/512  &  7.052E-02  &  0.35  &  1.872E-04  &  1.03    \\ 
1/1024  &  5.532E-02  &  0.35  &  9.157E-05  &  1.03    \\ 
1/2048  &  4.340E-02  &  0.35  &  4.482E-05  &  1.03    \\   \hline
Pred. &  & 0.35 &   & 0.70 \\
\hline
\end{tabular}
\caption{Example 2. Convergence rates for $\alpha = 1.4$ and $r = 0.3149$.}
\label{Table2}
\end{table}

\underline{Example 3}. \\
With $\alpha = 1.6$, $r = 1/2$, 
 \be
 u(x) \ = \ x^{\alpha/2} (1 - x) ^{\alpha/2} \, .
   \label{ppo15}
\ee
The corresponding right hand side is
\be
 f(x) \ = \ - \Gamma(1 + \alpha) \, \cos(\pi  \alpha  / 2) \, .
    \label{ppo16}
\ee
The numerical results are presented in Table \ref{Table3}.

\begin{table}[htb!]
\centering
\begin{tabular}{|c||c|c||c|c|} \hline
     $h$ & $ | u  -  u_{h} |_{H^{\alpha/2}(I)}$ & Cvg. rate &
	   $\| u  -  u_{h} \|_{L^{2}(I)}$ & Cvg. rate  \\
\hline
1/64  &  3.502E-02  &    &  6.559E-04  &     \\
1/128  &  2.461E-02  &  0.51  &  3.081E-04  &  1.09     \\
1/256  &  1.734E-02  &  0.50  &  1.479E-04  &  1.06     \\
1/512  &  1.224E-02  &  0.50  &  7.205E-05  &  1.04     \\
1/1024  &  8.651E-03  &  0.50  &  3.542E-05  &  1.02     \\
1/2048  &  6.115E-03  &  	0.50	  &  1.752E-05  &  1.02    \\  \hline
Pred. &  & 0.50 &   & 1.0 \\
\hline
\end{tabular}
\caption{Example 3. Convergence rates for $\alpha = 1.6$ and $r = 0.5$.}
\label{Table3}
\end{table}

\underline{Example 4}. \\
With $\alpha = 1.6$, $p = 0.9$, $q = \alpha - p$, 
$r  \, = \,  \sin(\pi (p + 1)) / ( \sin(\pi (p + 1)) - \sin(\pi (\alpha - p)) )$
\be
 u(x) \ = \ x^{p} (1 - x)^{q}  \, .
   \label{ppo19}
\ee
The corresponding right hand side is
\be
 f(x) \ = \ - (1 - r) \, \Gamma(1 + \alpha)  \frac{\sin(\pi \alpha)}{ \sin(\pi (\alpha - p))} \, .
    \label{ppo20}
\ee
The numerical results are presented in Table \ref{Table4}.
 
 \begin{table}[htb!]
\centering
\begin{tabular}{|c||c|c||c|c|} \hline
     $h$ & $ | u  -  u_{h} |_{H^{\alpha/2}(I)}$ & Cvg. rate &
	   $\| u  -  u_{h} \|_{L^{2}(I)}$ & Cvg. rate  \\
\hline
1/64  &  7.732E-02  &    &  7.083E-04  &   \\
1/128  &  5.827E-02  &  0.41  &  3.216E-04  &  1.14   \\
1/256  &  4.402E-02  &  0.40  &  1.485E-04  &  1.12   \\
1/512  &  3.331E-02  &  0.40  &  6.947E-05  &  1.10   \\
1/1024  &  2.522E-02  &  0.40  &  3.289E-05  &  1.08   \\
1/2048  &  1.910E-02  &  0.40  &  1.572E-05  &  1.06   \\  \hline
Pred. &  & 0.40 &   & 0.80 \\
\hline
\end{tabular}
\caption{Example 4. Convergence rates for $\alpha = 1.6$ and $r = 0.2764$.}
\label{Table4}
\end{table}
 
The numerical results are consistent with the theoretical predictions. Of particular note is that changing 
the convex combination of the adjoint operators in the definition of $\mcL_{r}^{\alpha}$, i.e., the factor $r$,
changes the regularity of the solution, and hence the convergence rate of the FEM approximation.

\setcounter{equation}{0}
\setcounter{figure}{0}
\setcounter{table}{0}
\setcounter{theorem}{0}
\setcounter{lemma}{0}
\setcounter{corollary}{0}
\section{Spectral type method for the solution of  
$\mcL_{r}^{\alpha} \, u \, = \, f $}
 \label{secSpec}
In this section we discuss a ``spectral type'' approximation method for the numerical solution of
$\mcL_{r}^{\alpha} \, u \, = \, f$. Central to the method is the following two results.
\begin{lemma} \label{lmaSpec1}
For $n \, = \, 0, 1, 2, \ldots $,
\begin{align}
 \mcL_{1/2}^{\alpha} \, x^{\alpha/2} (1 - x)^{\alpha/2} \, x^{n} &= \    \sum_{j = 0}^{n} a_{n , j} \, x^{j} 
 \, , \ \ \ \mbox{ where}  \nonumber \\
a_{n , j} \ =  \
  (-1)^{(n+1)}  \, \cos(\pi \, \alpha/2) \, \Gamma(1 + \alpha/2) 
   & \ \ \  \frac{ (-1)^{j} \Gamma(1 + \alpha + j)}%
{\Gamma(1 + \alpha/2 \, - \, n + j) \, \Gamma(1 + n - j) \, \Gamma(j + 1)}  \, .
 \label{spm1}
\end{align}
\end{lemma} 
\textbf{Proof}: \\
As $u(x) \, = \,  x^{\alpha/2} (1 - x)^{\alpha/2} \, x^{n}$ satisfies $u(0) = u(1) = 0$, then
\[
\mcL_{r}^{\alpha} \, u(x) \ = \  - D \, D \left( \bfD^{-(2 - \alpha)} \, + \, \bfD^{-(2 - \alpha) *}  \right) u(x) \, .
\]
Using Maple (see Figure \ref{mapleSPL1} in the appendix),
\be
\bfD^{-(2 - \alpha)} u(x) \ = \ 
\frac{\Gamma(1 + \alpha/2 \, + \, n)}{\Gamma(3 - \alpha/2 \, + \, n)} \, x^{n + 2 - \alpha/2} \, 
\mbox{}_{2}F_{1}(-\alpha/2 \, , \, 1 + \alpha/2 \, + \, n \, ; \, 3 - \alpha/2 \, + \, n \, ; \, x ) 
\, ,   \label{spm2}
\ee
and
\begin{align}
\bfD^{-(2 - \alpha) *} u(x) &=      
    \frac{\Gamma(-2 + \alpha/2 \, - \, n)}{\Gamma(-\alpha/2 \, - \, n)} \, x^{n + 2 - \alpha/2} \, 
\mbox{}_{2}F_{1}(-\alpha/2 \, , \, 1 + \alpha/2 \, + \, n \, ; \, 3 - \alpha/2 \, + \, n \, ; \, x )   \nonumber \\
&+ 
 (-1)^{n} \, \Gamma(1 + \alpha/2) \, 
 \sum_{k = 0}^{n + 2} \frac{(-1)^{k} \, \csc(\pi \alpha/2 \, + \, k \pi) \, \sin(\pi \alpha \, + \, k \pi)
 \Gamma(-1 + \alpha + k)}{\Gamma(-1 + \alpha/2 \, - n + k) \, \Gamma(3 + n - k) \, 
 \Gamma(k + 1)} \, x^{k} .
\label{spm3}
\end{align}
Using the identity
\be
  \Gamma(1 - z) \ = \ \frac{\pi}{\sin(\pi z)} \, \frac{1}{\Gamma(z)} \, ,
 \label{spm4}
\ee
with $z \ = \ 1 + \alpha/2 \, + n$, i.e., $1 - z \ = \ -\alpha/2 \, - n$,
\begin{align}
\Gamma(-\alpha/2 \, - n) &= \ \frac{\pi}{\sin(\pi (1 + \alpha/2 \, + n))} \, \frac{1}{\Gamma(1 + \alpha/2 \, + n)} 
 \ = \ \frac{\pi}{\sin(\pi \alpha/2) \, \cos(\pi (n + 1))} \, \frac{1}{\Gamma(1 + \alpha/2 \, + n)}  \nonumber  \\
 &= \frac{(-1)^{(n+1)} \, \pi}{\sin(\pi \alpha/2) \, \Gamma(1 + \alpha/2 \, + n)} \, .
 \label{spm5}
\end{align} 
Again using \eqref{spm4} with $z \ = \  3 - \alpha/2 \, + \, n$,
\begin{align}
\Gamma(-2 + \alpha/2 \, - n) &= \ \frac{\pi}{\sin(\pi (3 - \alpha/2 \, + n))} \, \frac{1}{\Gamma(3 - \alpha/2 \, + n)} 
 \ = \ \frac{\pi}{\sin(- \pi \alpha/2) \, \cos(\pi (n + 3))} \, \frac{1}{\Gamma(3 - \alpha/2 \, + n)}  \nonumber  \\
 &= \frac{(-1)^{(n+4)} \, \pi}{\sin(\pi \alpha/2) \, \Gamma(3 - \alpha/2 \, + n)} \, .
 \label{spm6}
\end{align} 
In view of \eqref{spm5} and \eqref{spm6}, we note that when adding  $\bfD^{-(2 - \alpha)} u(x)$ and
$\bfD^{-(2 - \alpha)*} u(x)$ the $x^{n + 2 - \alpha/2} \, \mbox{}_{2}F_{1}(\cdot)$ terms cancel.

Next, using standard trigonometric identities it is straightforward to show
\be
 \csc(\pi \alpha/2 \, + \, k \pi) \, \sin(\pi \alpha \, + \, k \pi) \ = \ 2 \cos(\pi \alpha/2) \, .
  \label{spm7}
\ee

Thus,
\begin{align*}
- D \, D \left( \frac{1}{2} \bfD^{-(2 - \alpha)} \, \right. &+  \left. \, \frac{1}{2} \bfD^{-(2 - \alpha) *}  \right) u(x)  \\
&=  \frac{1}{2} 
 (-1)^{(n + 1)} \, \Gamma(1 + \alpha/2) \, 
 \sum_{k = 2}^{n + 2} \frac{(-1)^{k} \, k \, (k - 1) \, 2 \, \cos(\pi \alpha/2) \, 
 \Gamma(-1 + \alpha + k)}{\Gamma(-1 + \alpha/2 \, - n + k) \, \Gamma(3 + n - k) \, 
 \Gamma(k + 1)} \, x^{(k - 2)} \, 
\end{align*}
which, after reindexing, yields \eqref{spm1}. \\
\mbox{ } \hfill \qed 

\begin{lemma} \label{lmaSpec2}
For $1 < \alpha < 2$, $0 \le \beta \le \alpha$, and $r$ satisfying
\be
     r \ = \ \frac{\sin(\pi \beta)}{\sin(\pi (\alpha - \beta)) \, + \, \sin(\pi \beta)} \, ,
\label{spm11}
\ee 
for $n \, = \, 0, 1, 2, \ldots $,
\begin{align}
 \mcL_{r}^{\alpha} \, x^{\beta} (1 - x)^{\alpha - \beta} \, x^{n} &= \    \sum_{j = 0}^{n} a_{n , j} \, x^{j} 
 \, , \ \ \ \mbox{ where}  \nonumber \\
a_{n , j} \ =  \ 
  (-1)^{(n+1)} (1 - r) \, \frac{\sin(\pi \, \alpha)}{\sin(\pi(\alpha - \beta))} \, \Gamma(1 + \alpha - \beta) 
  & \ \   \frac{ (-1)^{j} \, \Gamma(1 + \alpha + j)}%
{\Gamma(1 + \alpha - \beta \, - \, n + j) \, \Gamma(1 + n - j) \, \Gamma(j + 1)} \, .
 \label{spm12}
\end{align}
\end{lemma} 
\textbf{Proof}: \\       
With $u(x) \, = \,  x^{\beta} (1 - x)^{(\alpha - \beta)} \, x^{n}$ 
using Maple (see Figure \ref{mapleSPL2} in the appendix),
\be
\bfD^{-(2 - \alpha)} u(x) \ = \ 
\frac{\Gamma(1 + \beta \, + \, n)}{\Gamma(3 - \alpha + \beta \, + \, n)} \, x^{n + 2 - \alpha + \beta} \, 
\mbox{}_{2}F_{1}(1 + \beta + n \, , \, -\alpha  + \beta \, ; \, 3 - \alpha + \beta \, + \, n \, ; \, x ) 
\, ,   \label{spm13}
\ee
and
\begin{align}
\bfD^{-(2 - \alpha) *} u(x) &=      
    \frac{\Gamma(-2 + \alpha - \beta \, - \, n)}{\Gamma(-\beta \, - \, n)} \, x^{n + 2 - \alpha/2} \, 
\mbox{}_{2}F_{1}(1 + \beta + n \, , \, -\alpha  + \beta \, ; \, 3 - \alpha + \beta \, + \, n \, ; \, x )    \nonumber \\
&+ 
 (-1)^{n} \, \Gamma(1 + \alpha - \beta) \, 
 \sum_{k = 0}^{n + 2} \frac{(-1)^{k} \, \csc(\pi (\alpha - \beta) \, + \, k \pi) \, \sin(\pi \alpha \, + \, k \pi)
 \Gamma(-1 + \alpha + k)}{\Gamma(-1 + \alpha - \beta \, - n + k) \, \Gamma(3 + n - k) \, 
 \Gamma(k + 1)} \, x^{k} .
\label{spm14}
\end{align}
Using \eqref{spm5} with $\alpha/2 \longrightarrow  \beta$, and \eqref{spm6} with 
$\alpha/2 \longrightarrow (\alpha - \beta)$, we have that
\be
\frac{\Gamma(-2 + \alpha - \beta - n)}{\Gamma(-\beta - n)} 
\ = \ \frac{- \, \sin(\pi \beta) \, \Gamma(1 + \beta + n)}{\sin(\pi (\alpha - \beta)) \, \Gamma(3 - \alpha + \beta + n)}
\, .
   \label{spm15}
\ee   

The coefficient of $x^{n + 2 - \alpha + \beta} \mbox{}_{2}F_{1}(\cdot) $
in the linear combination $\left( r \, \bfD^{-(2 - \alpha)} \ + \ (1 - r) \, \bfD^{-(2 - \alpha) *}  \right) u(x)$ is:
\begin{align*}
r \frac{\Gamma(1 + \beta \, + \, n)}{\Gamma(3 - \alpha + \beta \, + \, n)}   &+ 
(1 - r) \frac{\Gamma(-2 + \alpha - \beta \, - \, n)}{\Gamma(-\beta \, - \, n)}  \\
&= \frac{\Gamma(1 + \beta \, + \, n)}{\Gamma(3 - \alpha + \beta \, + \, n)}
\left(r \ + \ (1 - r) \frac{- \, \sin(\pi \beta)}{\sin(\pi (\alpha - \beta))} \right) \ \ \mbox{(using \eqref{spm15})} \\
&= 0 \, ,
\end{align*}
provided $r$ is given by \eqref{spm11}.

Using standard trigonometric identities it is straightforward to show
\[
 \csc(\pi (\alpha - \beta) \, + \, k \pi) \, \sin(\pi \alpha \, + \, k \pi) \ = \ 
 \frac{\sin(\pi \alpha)}{\sin(\pi (\alpha - \beta))} \, .
\]

Thus,
\begin{align*}
- D \, D \left( r \bfD^{-(2 - \alpha)} \, +   \, (1 - r) \bfD^{-(2 - \alpha) *}  \right) u(x)  & \\
=  \
 (-1)^{(n + 1)} (1 - r) \, \Gamma(1 + \alpha - \beta) \, & 
 \sum_{k = 2}^{n + 2} \frac{(-1)^{k} \, k \, (k - 1) \, \frac{\sin(\pi \alpha)}{\sin(\pi (\alpha - \beta))} \, 
 \Gamma(-1 + \alpha + k)}{\Gamma(-1 + \alpha - \beta \, - n + k) \, \Gamma(3 + n - k) \, 
 \Gamma(k + 1)} \, x^{(k - 2)} \, 
\end{align*}
which, after reindexing, yields \eqref{spm12}. \\
\mbox{ } \hfill \qed 
       
Jacobi polynomial play a key role in the approximation schemes. We briefly review their definition and
properties central to the method \cite{abr641, sze751}. 

\textbf{Usual Jacobi Polynomials, $P_{n}^{(\alpha , \beta)}(x)$, on $(-1 \, , \, 1)$}.   \\    
\underline{Definition}: $ P_{n}^{(\alpha , \beta)}(x) \ := \ 
\sum_{m = 0}^{n} \, p_{n , m} \, (x - 1)^{(n - m)} (x + 1)^{m}$, where
\be
       p_{n , m} \ := \ \frac{1}{2^{n}} \, \left( \begin{array}{c}
                                                              n + \alpha \\
                                                              m  \end{array} \right) \,
                                                    \left( \begin{array}{c}
                                                              n + \beta \\
                                                              n - m  \end{array} \right) \, .
\label{spm21}
\ee
\underline{Orthogonality}:    
\begin{align}
 & \int_{-1}^{1} (1 - x)^{\alpha} (1 + x)^{\beta} \, P_{j}^{(\alpha , \beta)}(x) \, P_{k}^{(\alpha , \beta)}(x)  \, dx 
 \ = \
   \left\{ \begin{array}{ll} 
   0 , & k \ne j \, , \\
   |\| P_{j}^{(\alpha , \beta)} |\|^{2}
   \, , & k = j  
   \, . \end{array} \right.   \nonumber \\
& \quad \quad \mbox{where } \  \ |\| P_{j}^{(\alpha , \beta)} |\| \ = \
 \left( \frac{2^{(\alpha + \beta + 1)}}{(2j \, + \, \alpha \, + \, \beta \, + 1)} 
   \frac{\Gamma(j + \alpha + 1) \, \Gamma(j + \beta + 1)}{\Gamma(j + 1) \, \Gamma(j + \alpha + \beta + 1)}
   \right)^{1/2} \, .
  \label{spm22}
\end{align}                                                    
       
 \textbf{Jacobi Polynomials, $G_{n}(p , q , x)$, on $(0 \, , \, 1)$}.   \\    
\underline{Definition}: $G_{n}(p , q , x) \ := \ 
\sum_{j = 0}^{n} \, g_{n , j} \, x^{j}$, where
\be
       g_{n , j} \ := \ (-1)^{(n-j)} \frac{\Gamma(q + n)}{\Gamma(p + 2n)} \,
         \frac{\Gamma(n + 1)}{\Gamma(j + 1) \, \Gamma(n - j + 1)} \, 
         \frac{\Gamma(p + n + j)}{\Gamma(q + j)} \, .
\label{spm31}
\ee
\underline{Orthogonality}:    
\begin{align}
 & \int_{0}^{1} x^{(q - 1)} (1 - x)^{(p - q)} \, G_{j}(p , q , x) \, G_{k}(p , q , x)  \, dx \ = \
   \left\{ \begin{array}{ll} 
   0 , & k \ne j \, , \\
   |\| G_{j}^{(p , q)} |\|^{2}
   \, , & k = j  
   \, . \end{array} \right.   \nonumber \\
& \quad \quad \mbox{where } \  \ |\| G_{n}^{(p , q)} |\| \ = \
 \left( \frac{\Gamma(n + 1) \, \Gamma(n + q) \, \Gamma(n + p) \, \Gamma(n + p - q + 1)}%
 {(2n + p) \, \Gamma^{2}(2n + p)} 
   \right)^{1/2} \, .
  \label{spm32}
\end{align}      

Note that $G_{n}(p , q , x) \ = \ \frac{\Gamma(n + 1) \, \Gamma(n + p)}{\Gamma(2n + p)} \,
P_{n}^{(p - q \,  , \, q - 1)}(2x - 1)$.

\textbf{The weighted $L^{2}(0 , 1)$ spaces, $L_{\rho}^{2}(0 , 1)$}. \\
The weighted $L^{2}(0 , 1)$ spaces are convenient for analyzing the convergence of the spectral type methods 
presented below. For $\rho(x) > 0, \ x \in (0 , 1)$, let 
\[
L_{\rho}^{2}(0 , 1) \, := \, \{ f(x) \, : \, \int_{0}^{1} \rho(x) \, f(x)^{2} \, dx \ < \ \infty \} \, .
\]
Associated with $L_{\rho}^{2}(0 , 1)$ is the inner product, $\langle \cdot , \cdot \rangle_{\rho}$, and
norm, $\| \cdot \|_{\rho}$, defined by
\begin{align*}
\langle f \,  , \, g \rangle_{\rho} &:= \ \int_{0}^{1} \rho(x) \, f(x) \, g(x) \, dx \, , \quad \mbox{and} \\
 \| f \|_{\rho} &:= \ \left( \langle f \,  , \, f \rangle_{\rho} \right)^{1/2} \, .
\end{align*}

\subsection{Spectral type method approximation to $\mcL_{1/2}^{\alpha} u \ = \ f$}
\label{ssecsp1}
In this section we discuss the approximation of $\mcL_{1/2}^{\alpha} u \ = \ f$, subject to $u(0) = u(1) = 0$, using 
Jacobi polynomials on $(0 , 1)$. For $1 < \alpha < 2$ (fixed), for convenience of notation, we let
$G_{n}(x) := G_{n}(1 + \alpha \, , \, 1 + \alpha/2 \, , \, x)$, and let $\mcP_{n}(x)$ denote the vector space
of polynomials of degree $\le n$. Let
\begin{align*}
   \omega(x) &= \ x^{\alpha/2} (1 - x)^{\alpha/2},  \\
\mbox{and }  \quad \lambda_{n} &= \ - \cos(\pi \, \alpha / 2) 
 \, \frac{\Gamma(n + 1 + \alpha)}{\Gamma(n + 1)} \, . 
\end{align*}
We have that 
\[
\| G_{n} \|_{\omega}^{2}  \ = \ \int_{0}^{1} x^{\alpha/2} (1 - x)^{\alpha/2} \, G_{n}(x) \, G_{n}(x) \, dx \ = \
  |\| G_{n}^{(1 + \alpha \, , \, 1 + \alpha/2 )} |\|^{2} \, .
\]  

Additionally, $\{ G_{j}(x) \}_{j = 0}^{\infty}$ is an orthogonal basis for $L^{2}_{\omega}(0 , 1)$.

Using Stirling's formula we have that
\be
\lim_{n \rightarrow \infty} \, \frac{\Gamma(n + \mu)}{\Gamma(n) \, n^{\mu}}
\ = \ 1 \, , \mbox{ for } \mu \in \real.  
 \label{eqStrf}
\ee 
Thus $\lambda_{n} > 0$ for all $n \, = \, 0, 1, 2, \ldots $, and as $n \rightarrow \infty$
$\lambda_{n} \sim - \cos(\pi \, \alpha / 2) \, (n + 1)^{\alpha}$.
 
\textbf{Remark}: Note that $f \in L^{2}_{\omega}(0 , 1)$ may be expressed as $f(x) \ = \ \sum_{i = 0}^{\infty} 
\frac{f_{i}}{ \| G_{i} \|_{\omega}^{2} } \, G_{i}(x)$, where $f_{i}$ is given by
\be
f_{i} \, := \, \int_{0}^{1} \, \omega(x) \, G_{i}(x) \, f(x) \, dx \, .
\label{deffi}
\ee

We begin with the following important extension of Lemma \ref{lmaSpec1}.
\begin{lemma} \label{lmaGnO}
For $n \, = \, 0, 1, 2, \ldots $,
\be
 \mcL_{1/2}^{\alpha} \,\omega(x) \, G_{n}(x) \ = \  \lambda_{n} \, G_{n}(x) \, .
 \label{eeGO}
\ee
\end{lemma}

\textbf{Proof}: We have that, up to a constant multiplier, $G_{n}(x)$ is characterized by
\begin{align*}
0 &= \left( \omega(x) \, G_{n}(x) \ , \ p(x) \right) \, , \ \forall \, p(x) \,  \in \, \mcP_{n-1}(x) \, , \\
&= \int_{0}^{1} \omega(x) \, G_{n}(x) \, p(x) \, dx \, .
\end{align*}
 
Now, given $p(x) \in  \mcP_{n-1}(x)$, from Lemma \ref{lmaSpec1}, there exists $\hat{p}(x) \in  \mcP_{n-1}(x)$
satisfying
\[
\hat{p}(x) \ = \ \mcL_{1/2}^{\alpha} \, \omega(x) \, p(x) \, .
\]
Then, noting that $\mcL_{1/2}^{\alpha}$ is self adjoint, 
\begin{align*}
\left( \omega(x) \, \mcL_{1/2}^{\alpha} \, \omega(x) \, G_{n}(x) \ , \ p(x) \right)
 &=  \left( \mcL_{1/2}^{\alpha} \, \omega(x) \, G_{n}(x) \ , \  \omega(x) \, p(x) \right) \\
 &=  \left( \omega(x) \, G_{n}(x) \ , \  \mcL_{1/2}^{\alpha}  \, \omega(x) \, p(x) \right) \\
 &=  \left( \omega(x) \, G_{n}(x) \ , \  \hat{p}(x) \right) \\
 &= 0 \, ,
 \end{align*}
 which implies that $\mcL_{1/2}^{\alpha} \, \omega(x) \, G_{n}(x) \ = \ C \, G_{n}(x)$, for $C \in \real$.
 
 As the coefficient of $x^{n}$ in $G_{n}(x)$ is $1$, then from Lemma \ref{lmaSpec1},
 \[
 C \ = \ a_{n , n} \ = \ - \cos(\pi \, \alpha / 2) 
 \frac{\Gamma(1 + \alpha + n)}{\Gamma(n + 1)} \ = \ \lambda_{n} \, .
 \]
\mbox{ } \hfill \qed

With $f_{i}$ defined in \eqref{deffi}, let 
\be
u_{N}(x) \ := \  \omega(x) \, \sum_{j = 0}^{N} c_{j} \, G_{j}(x) \, , \mbox{ where} \ 
c_{j} \ = \  \frac{1}{\lambda_{j} \, \| G_{j}  \|_{\omega}^{2} } \, f_{j}    \, .
\label{defuN}
\ee

\begin{theorem}  \label{thmusps}
Let $f (x) \in L^{2}_{\omega}(0 , 1)$ and $u_{N}(x)$ be as defined in \eqref{defuN}. Then,
$u(x) \ := \ \lim_{N \rightarrow \infty} u_{N}(x) \ = \ \omega(x) \, \sum_{j = 0}^{\infty}  c_{j} \, G_{j}(x) \, 
\in L^{2}_{\omega^{-1}}(0 , 1)$.
In addition, $\mcL_{1/2}^{\alpha} u(x) \ = \ f(x)$.
\end{theorem}
\textbf{Proof}: 
For $f_{N}(x) \ = \ \sum_{i = 0}^{N} 
\frac{f_{i}}{ \| G_{i} \|_{\omega}^{2} } \, G_{i}(x)$, we have that $f(x) \ = \ \lim_{N \rightarrow \infty} f_{N}(x)$, and 
$\left\{ f_{N}(x) \right\}_{N = 0}^{\infty}$ is a Cauchy sequence in $L^{2}_{\omega}(0 , 1)$. 
A straightforward calculation shows that $u_{N}(x) \, \in L^{2}_{\omega^{-1}}(0 , 1)$.
Then, (without loss of 
generality, assume $M > N$)
\begin{align*}
\| u_{N}(x) \ - \ u_{M}(x) \|_{\omega^{-1}}^{2} 
&= \ \left( \omega^{-1}(x) \ \omega(x) \, \sum_{j = N + 1}^{M} c_{j} \, G_{j}(x) \ , \ 
   \omega(x) \, \sum_{j = N + 1}^{M} c_{j} \, G_{j}(x) \right)     \\
&= \
\left(  \omega(x) \, \sum_{j = N + 1}^{M} \frac{f_{j}}{\lambda_{j} \, \| G_{j}  \|_{\omega}^{2} }  \, G_{j}(x) \ , \ 
   \sum_{j = N + 1}^{M} \frac{f_{j}}{\lambda_{j} \, \| G_{j}  \|_{\omega}^{2} }  \, G_{j}(x) \right)  \\
&= \
 \sum_{j = N + 1}^{M} \frac{f_{j}^{2}}{\lambda_{j}^{2} \, \| G_{j}  \|_{\omega}^{2} }   \\
&\le \ C \,
\left(  \omega(x) \, \sum_{j = N + 1}^{M} \frac{f_{j}}{ \| G_{j}  \|_{\omega}^{2} }  \, G_{j}(x) \ , \ 
   \sum_{j = N + 1}^{M} \frac{f_{j}}{ \| G_{j}  \|_{\omega}^{2} }  \, G_{j}(x) \right)  \\
& \quad \quad \mbox{ (using $\lambda_{j}$'s are bounded away from zero) }  \\
&= \ C \, \| f_{N}(x) \ - \ f_{M}(x) \|_{\omega}^{2} \, .
\end{align*}
Hence $\{ u_{N}(x) \}_{N = 0}^{\infty}$ is a Cauchy sequence in $L^{2}_{\omega^{-1}}(0 , 1)$. As 
$L^{2}_{\omega^{-1}}(0 , 1)$ is closed,   \linebreak[4]
$u(x) \ := \ \lim_{N \rightarrow \infty} u_{N}(x) \, \in \, L^{2}_{\omega^{-1}}(0 , 1)$.

Next, as $f_{N}(x) \rightarrow f(x)$ in $L^{2}_{\omega}(0 , 1)$, given $\eps > 0$ there exists $\tilde{N}$ 
such that for $N > \tilde{N}$, $\| f(x)  \ - \ f_{N}(x) \|_{\omega} \, < \, \eps$. Then, for $N > \tilde{N}$,
using Lemma \ref{lmaGnO}
\begin{align*}
\| f(x) \ - \ \mcL^{\alpha}_{1/2} u_{N}(x) \|_{\omega} &= \| f(x) \ - \  
\mcL^{\alpha}_{1/2} \left(  \omega(x) \, \sum_{j = 0}^{N}  \frac{f_{j}}{\lambda_{j} \, \| G_{j}  \|_{\omega}^{2} }
 \, G_{j}(x)  \right) \|_{\omega} \\
&= \ 
\| f(x) \ - \  
\sum_{j = 0}^{N}  \frac{f_{j}}{ \| G_{j}  \|_{\omega}^{2} }
 \, G_{j}(x)  \|_{\omega} \\
&= \ 
\| f(x) \ - \  f_{N}(x)  \|_{\omega} \ < \ \eps \, .
\end{align*}
Hence, $f(x) \ = \ \mcL^{\alpha}_{1/2} u(x)$. \\
\mbox{ } \hfill \qed

Using Lemma \ref{lmaGnO} and $B(\cdot , \cdot)$ defined in \eqref{treg1p5}, we have the following 
connection between the spectral type approximation and the Galerkin approximation.
\begin{lemma} \label{solspcse}
Let $f(x) \in L^{2}_{\omega}(0 , 1)$ and 
 $X_{N} \, := \, \left\{ h(x) \, : \, h(x) \ = \ \omega(x) \, \tilde{h}(x) \, , 
\ \tilde{h}(x) \in \mcP_{N}(x) \right\}$. Then, $u_{N}(x) \in X_{N}$ satisfying
$B(u_{N} , v) \ = \ \langle f , v \rangle \ \forall v \in X_{N}$, is given by \eqref{defuN}.
\end{lemma}
\textbf{Remark}: Note that the factor $\omega(x) \ = \ x^{\alpha/2} (1 - x)^{\alpha/2}$ explicitly incorporates the
endpoint singular behavior of the operator into the approximation space $X_{N}$.

\textbf{Proof}: The set $\{ G_{i}(x) \}_{i = 0}^{N}$ forms a basis for $\mcP_{N}(x)$. With 
$u_{N}(x) \ = \ \omega(x) \, \sum_{j = 0}^{N} c_{j} \, G_{j}(x)$, $v(x) \ = \ \omega(x) \, G_{i}(x)$, from
\eqref{deffi} (using \eqref{eeGO}) we obtain
\begin{align*}
f_{i} &=  \ \left( \omega(x) \, G_{i}(x) \ , \ f(x) \right) \ =  \ \langle f \, , \, v \rangle \\
&= \ B(u_{N} \, , \, v) \ = \ \left( \mcL_{1/2}^{\alpha}  \, \omega(x) \, \sum_{j = 0}^{N} c_{j} \, G_{j}(x) \ , \
 \omega(x) \, G_{i}(x) \right)   \\
 &= \ \left(  \sum_{j = 0}^{N} c_{j} \, \lambda_{j} \, G_{j}(x) \ , \
 \omega(x) \, G_{i}(x) \right)   \\
 &= \ c_{i} \, \lambda_{i} \, \| G_{i} \|_{\omega}^{2} \, .
\end{align*}

Hence, $c_{i} \ = \ f_{i} / ( \lambda_{i} \,  \| G_{i} \|_{\omega}^{2} ) $. \\
\mbox{ } \hfill  \qed

\subsubsection{A priori error estimate for $u \, - \, u_{N}$}
\label{sssecAP1}
We have the following statement for the error between $u \, - \, u_{N}$.
\begin{theorem}   \label{apespcse}
For $f(x) \in L^{2}_{\omega}(0 , 1)$ and $u_{N}(x)$ given by  \eqref{defuN}, there exists $C > 0$ such that
\be
 \| u \, - \, u_{N} \|_{\omega^{-1}} \ \le \   \frac{1}{ \lambda_{N+1} } \, \| f \|_{\omega} \ \le \
 \ C \, (N + 2)^{- \alpha}    \, \| f \|_{\omega} . 
  \label{apert1}
\ee 
\label{thmAP1}
\end{theorem} 
\textbf{Proof}: With the definition of the $ \| \cdot \|_{\omega^{-1}} $ norm, and \eqref{eqStrf}
\begin{align*}
 \| u \, - \, u_{N} \|_{\omega^{-1}}^{2} &=
 \int_{0}^{1} \omega^{-1}(x) \left( \omega(x) \, \sum_{i = 0}^{\infty}  \frac{ f_{i} }{ ( \lambda_{i} \,  \| G_{i} \|_{\omega}^{2} )}
  \, G_{i}(x) 
 \ - \ \omega(x) \, \sum_{i = 0}^{N}  \frac{ f_{i} }{ ( \lambda_{i} \,  \| G_{i} \|_{\omega}^{2} ) } \, G_{i}(x) \right)^{2} \, dx  \\
&\le \ \max_{N + 1 \, \le \, i} \left( \frac{1}{  \lambda_{i}^{2}  } \right) \ 
  \int_{0}^{1} \, \omega(x) \, \sum_{i \, = \, N + 1}^{\infty} \left( \frac{ f_{i} }{\| G_{i} \|_{\omega}^{2} } \, G_{i}(x) 
  \right)^{2} \, dx  \\
&\le \  \frac{1}{  \lambda_{N + 1}^{2}  }  \ 
  \int_{0}^{1} \, \omega(x) \, \sum_{i \, = \, 0}^{\infty} \left( \frac{ f_{i} }{\| G_{i} \|_{\omega}^{2} } \, G_{i}(x) 
  \right)^{2} \, dx  \\
&\le \  \left( \frac{1}{ - \cos(\pi \, \alpha / 2)} \, \frac{\Gamma(N + 2)}{\Gamma(N + 2 + \alpha)} \right)^{2} \ 
  \int_{0}^{1} \, \omega(x) \, \sum_{i \, = \, 0}^{\infty} \left( \frac{ f_{i} }{\| G_{i} \|_{\omega}^{2} } \, G_{i}(x) 
  \right)^{2} \, dx  \\
&= \  \left( \frac{1}{ - \cos(\pi \, \alpha / 2)} \, \frac{\Gamma(N + 2)}{\Gamma(N + 2 + \alpha)} \right)^{2} \ 
  \int_{0}^{1} \, \omega(x) \, f(x)^{2} \, dx  \\
&\le \ \left( \frac{1}{ - \cos(\pi \, \alpha / 2)} \, \frac{\Gamma(N + 2)}{\Gamma(N + 2 + \alpha)}  \right)^{2} \ 
  \| f \|_{\omega}^{2}  \\
&\le \ C \, (N + 2)^{-2 \, \alpha}   \ 
  \| f \|_{\omega}^{2} \, .
\end{align*}
\mbox{ } \hfill \qed

\begin{corollary}
For $f(x) \in L^{2}_{\omega}(0 , 1)$ and $u_{N}(x)$ given by  \eqref{defuN}, there exists $C > 0$ such that
\begin{align}
 \| u \, - \, u_{N} \|_{H^{\alpha/2}(0 , 1)} &\le \   \frac{1}{\sqrt{ \lambda_{N+1} } } \, \| f \|_{\omega} \ \le \
 \ C \, (N + 2)^{- \alpha / 2}    \, \| f \|_{\omega} \, ,   \mbox{ and } \label{apert2}   \\
 \| u \, - \, u_{N} \| &\le \   \frac{1}{\lambda_{N+1} } \, \| f \|_{\omega} \ \le \
 \ C \, (N + 2)^{- \alpha }    \, \| f \|_{\omega} \, . \label{apert3}
\end{align} 
\end{corollary} 
\label{corAP2}
\textbf{Proof}: With $\langle \cdot \, , \, \cdot \rangle$ denoting the $L^{2}$-duality pairing,
and using the coercivity of $B( \cdot , \cdot)$ (see \cite{erv061}), there exists $C_{0} > 0$ such that
\begin{align*}
C_{0} \,  \| u \, - \, u_{N} \|_{H^{\alpha/2}(0 , 1)}^{2} &\le  \ B(u - u_{N} \, , \, u - u_{N})  \nonumber \\
 &= \ \langle \mcL_{1/2}^{\alpha} \, \omega(x) \, \sum_{j \, = \, N + 1}^{\infty} \frac{ f_{j} }{\lambda_{j} \, 
 \| G_{j} \|_{\omega}^{2}}
 \, G_{j}(x) \ , \ \omega(x) \, \sum_{j \, = \, N + 1}^{\infty} \frac{ f_{j} }{\lambda_{j} \, \| G_{j} \|_{\omega}^{2}} 
 \, G_{j}(x) \rangle  
 \nonumber \\
&= \ \left(  \sum_{j \, = \, N + 1}^{\infty} \frac{ f_{j} }{\lambda_{j} \, \| G_{j} \|_{\omega}^{2}}
\, \lambda_{j} \, G_{j}(x) \ , \ \omega(x) \, \sum_{j \, = \, N + 1}^{\infty} \frac{ f_{j} }{\lambda_{j} \, \| G_{j} \|_{\omega}^{2}} \, G_{j}(x) \right)  
 \nonumber \\
&= \ \sum_{j \, = \, N + 1}^{\infty}  \,  \frac{ f_{j}^{2} }{\lambda_{j} \, \| G_{j} \|_{\omega}^{2}}  \nonumber \\
&\le \  \max_{N + 1 \, \le \, j} \left( \frac{1}{ \lambda_{j} } \right) \, 
    \sum_{j \, = \, 0}^{\infty}  \,  \frac{ f_{j}^{2} }{ \| G_{j} \|_{\omega}^{2}}  \nonumber \\
&\le \ \frac{1}{ \lambda_{N + 1} }   \, \| f \|_{\omega}^{2}   \nonumber \\
&\le \ C \, (N + 2)^{- \alpha}   \ 
  \| f \|_{\omega}^{2} \, ,  
\end{align*}
where in the last step we use the bound for $\lambda_{N + 1}$ obtained in the proof of Theorem \ref{thmAP1}.

The bound \eqref{apert3} for the $L^{2}$ error in the approximation follows from \eqref{apert1} and the observation that, as
$\omega(x) \, = \, x^{\alpha/2} (1 - x)^{\alpha/2} \, < 1 $, for $0 < x < 1$, 
$\| u \, - \, u_{N} \| \ \le \ \| u \, - \, u_{N} \|_{\omega^{-1}}$.  \\
\mbox{ } \hfill \qed

            
\subsection{Spectral type method approximation to $\mcL_{r}^{\alpha} u \ = \ f$}
\label{ssecsp2}
For the general case  $\mcL_{r}^{\alpha} u \ = \ f$, $r \ne 1/2$, the operator $\mcL_{r}^{\alpha} \cdot $ is not
symmetric. Hence the singular behavior of the adjoint problem $ \left( \mcL_{r}^{\alpha} \right)^{*} \cdot \ = \
\mcL_{1 - r}^{\alpha} \cdot$ does not match that of $\mcL_{r}^{\alpha} \cdot$. In order to conveniently
present the approximation method and its properties, in this section we use the following notation.

For $1 < \alpha < 2$ and $r$ given, and $\beta$ determined by \eqref{spm11},
\begin{align}
& \mcL_{r}^{\alpha} u  \ = \ r \, D^{\alpha} u \ + \ (1 - r) \, D^{\alpha *} u 
 & \quad \quad &  
 \mcL_{r}^{\alpha *} u \ = \ r \, D^{\alpha *} u \ + \ (1 - r) \, D^{\alpha} u   \nonumber  \\
& \omega(x)  \ = \ x^{\beta} \, (1 - x)^{\alpha - \beta} 
 & \quad \quad &  
\omega^{*}(x) \ = \ x^{\alpha - \beta} \, (1 - x)^{\beta}   \label{wwstar}  \\
& \mcG_{n}(x)  \ = \ G_{n}(\alpha + 1 \, , \, \beta + 1 \, , \, x)  
 & \quad \quad &  
\mcG_{n}^{*}(x) \ = \ G_{n}(\alpha + 1 \, , \, \alpha - \beta + 1 \, , \, x)  \nonumber  \\
& \lambda_{n}  \ = \ -(1 - r) \frac{\sin( \pi \, \alpha)}{\sin( \pi (\alpha - \beta) )} \,
\frac{\Gamma(n + 1 + \alpha)}{\Gamma(n + 1)}
 & \quad \quad &  
\lambda_{n}^{*} \ = \ -r \frac{\sin( \pi \, \alpha)}{\sin( \pi (\alpha - \beta) )} \,
\frac{\Gamma(n + 1 + \alpha)}{\Gamma(n + 1)}            \nonumber  
\end{align}
From \eqref{spm32} we have the following orthogonality properties
\[
\int_{0}^{1} \omega(x) \, \mcG_{j}(x) \, \mcG_{k}(x) \, dx \ = \ 0 \, , \ \
k \ne j \, , \quad \quad
\int_{0}^{1} \omega^{*}(x) \, \mcG_{j}^{*}(x) \, \mcG_{k}^{*}(x) \, dx \ = \ 0 \, , \ \
k \ne j \, ,
\]
and,
\begin{align}
\| \mcG_{n} \|_{\omega}^{2} &= \ \| G_{n}(\alpha + 1 \, , \, \beta + 1 \, , \, x) \|_{\omega}^{2} \nonumber \\
&= \ 
\frac{\Gamma(n + 1) \, \Gamma(n + \alpha + 1) \, \Gamma(n + \beta + 1) \, \Gamma(n + \alpha - \beta + 1)}%
{ (2 n \, + \, \alpha \, + \, 1) \ \Gamma^{2}(2 n \, + \, \alpha \, + \, 1) }  \nonumber \\
&= \ \| G_{n}(\alpha + 1 \, , \, \alpha - \beta + 1 \, , \, x) \|_{\omega^{*}}^{2}  \nonumber \\
&= \  \| \mcG_{n}^{*} \|_{\omega^{*}}^{2} \, .  \label{dtgh1}
\end{align}

Corresponding to Lemma \ref{lmaGnO} we have the following.
\begin{lemma} \label{lmaGngO}
For $n \, = \, 0, 1, 2, \ldots $,
\begin{align}
 \mcL_{r}^{\alpha} \,\omega(x) \, \mcG_{n}(x) &= \  \lambda_{n} \, \mcG_{n}^{*}(x) \, ,
 \label{eeGgO}  \\
 \mcL_{1 - r}^{\alpha} \,\omega^{*}(x) \, \mcG_{n}^{*}(x) &= \  \lambda_{n}^{*} \, \mcG_{n}(x) \, .
 \label{eeGgO1}  
\end{align}
\end{lemma}

\textbf{Proof}: Up to a multiplicative constant, $\mcG_{n}(x)$ and $\mcG_{n}^{*}(x)$ are,
respectively, determined by  \linebreak[4]
$\left( \mcG_{n}(x) \, , \, p(x) \right)_{\omega} \ = \ 0 $ and
$\left( \mcG_{n}^{*}(x) \, , \, p(x) \right)_{\omega^{*}} \ = \ 0 $, for all $p(x) \in \mcP_{n-1}(x)$.

Let $p(x) \in \mcP_{n-1}(x)$. Then, from Lemma \ref{lmaSpec2} there exists $\hat{p}(x) \in P_{n-1}(x)$
such that $\mcL_{1 - r}^{\alpha} \omega^{*}(x) \, p(x) \ = \ \hat{p}(x)$. Then,
\begin{align*}
\left( \mcL_{r}^{\alpha} \left( \omega(x) \, \mcG_{n}(x) \right) \ , \ p(x) \right)_{\omega^{*}}
&= \ \int_{0}^{1} \omega^{*}(x) \, \mcL_{r}^{\alpha} \left( \omega(x) \, \mcG_{n}(x) \right) \, p(x) \, dx  \\
&= \ 
  \int_{0}^{1} \mcL_{r}^{\alpha} \left( \omega(x) \, \mcG_{n}(x) \right) \, \omega^{*}(x) \, p(x) \, dx  \\
&= \ 
  \int_{0}^{1}   \omega(x) \, \mcG_{n}(x)  \, \mcL_{1 - r}^{\alpha} \left( \omega^{*}(x) \, p(x) \right) \, dx  \\
&= \ 
  \int_{0}^{1}   \omega(x) \, \mcG_{n}(x)  \, \hat{p}(x) \, dx  \\
&= \ 0 \, .
\end{align*}
Hence, $\mcL_{r}^{\alpha} \, \omega(x) \, \mcG_{n}(x) \ = \ C \, \mcG_{n}^{*}(x)$, for $C \in \real$.
 
 As the coefficient of $x^{n}$ in $\mcG_{n}(x)$ and $\mcG_{n}^{*}(x)$ is $1$, then from Lemma \ref{lmaSpec2},
 \[
 C \ = \ -(1 - r) \frac{\sin( \pi \, \alpha)}{\sin( \pi (\alpha - \beta) )} \,
\frac{\Gamma(n + 1 + \alpha)}{\Gamma(n + 1)} \ = \ \lambda_{n} \, .
 \]
 
 An analogous argument to the above establishes \eqref{eeGgO1}. \\
\mbox{ } \hfill \qed

\textbf{Remark}: Note that $f(x) \in L_{\omega^{*}}^{2}(0 , 1)$ may be expressed as 
$f(x) \ = \ \sum_{i = 0}^{\infty} \frac{f_{i}^{*}}{ \| \mcG_{i}^{*} \|_{\omega^{*}}^{2}} \, \mcG_{i}^{*}(x)$, 
where $f_{i}^{*}$ is given by
\be
f_{i}^{*} \, := \, \int_{0}^{1} \,  \omega^{*}(x) \, f(x) \, \mcG_{i}^{*}(x)  \, dx \, .
\label{deffis}
\ee

With $f_{i}^{*}$ defined in \eqref{deffis}, let 
\be
u_{N}(x)  \ = \ \omega(x) \, \sum_{i = 0}^{N} c_{i} \, \mcG_{i}(x) \, , \mbox{ where} \ 
c_{i} \ = \ \frac{1}{ \lambda_{i} \, \| \mcG_{i}^{*} \|_{\omega^{*}}^{2}} f_{i}^{*} \, .
\label{defuNs}
\ee

\begin{theorem}   \label{thmuspg}
Let $f (x) \in L^{2}_{\omega^{*}}(0 , 1)$ and $u_{N}(x)$ be as defined in \eqref{defuNs}. Then,
$u(x) \ := \ \lim_{N \rightarrow \infty} u_{N}(x) \ = \ \omega(x) \, \sum_{j = 0}^{\infty}  c_{j} \, \mcG_{j}(x) \, 
\in L^{2}_{\omega^{-1}}(0 , 1)$.
In addition, $\mcL_{r}^{\alpha} u(x) \ = \ f(x)$.
\end{theorem}
\textbf{Proof}: Using \eqref{dtgh1}, the proof follows in a similar manner to that for Theorem \ref{thmusps}. \\
\mbox{ } \hfill \qed

\subsubsection{Invertibility of $\mcL_{r}^{\alpha} \cdot$ on $L^{2}(0 , 1)$}
\label{sssecL2inv}
We return to the question eluded to by Lemmas \ref{lmaK1} and \ref{lmamap2} in Section \ref{sec_reg}, 
namely the invertibility of $\mcL_{r}^{\alpha} \cdot$ on $L^{2}(0 , 1)$. Theorem \ref{thmuspg}, together with
\eqref{defuNs} and \eqref{deffis} (see also Theorem \ref{thmusps}, together with \eqref{defuN} and
\eqref{deffi}) gives an explicit inverse for $\mcL^{\alpha}_{r} \cdot$ on $L^{2}_{\omega^{*}}(0 , 1) \supset
L^{2}(0 , 1)$. Hence we have the following.
\begin{corollary} \label{invofL2}
For $1 < \alpha < 2$, $0 < r < 1$, $\beta$ chosen such that \eqref{spm11} is satisfied, $\omega$ and 
$\omega^{*}$ as in \eqref{wwstar}, given $f \in L^{2}(0 , 1)$ there exists a unique solution 
$u \in L^{2}_{\omega^{-1}}(0 , 1)$ such that $\mcL^{\alpha}_{r} u \, = \, f$ and $u(0) \, = \, u(1) \, = \, 0$.
(For a solution to the nonhomogeneous boundary condition problem: $\mcL^{\alpha}_{r} u_{nh} \ = \ f$
subject to $u_{nh}(0) = A$, $u_{nh}(1) = B$, the homogeneous boundary condition 
for $u$ is combined with a
suitable function chosen from the kernel of $\mcL^{\alpha}_{r} \cdot$ (see Corollary \ref{cor1r})).
\end{corollary}
\mbox{} \hfill \qed

\subsubsection{A priori error estimate for $u \, - \, u_{N}$}
\label{sssecAP2}
We have the following statement for the error between $u \, - \, u_{N}$.
\begin{theorem}
For $f(x) \in L^{2}_{\omega^{*}}(0 , 1)$ and $u_{N}(x)$ given by  \eqref{defuNs}, there exists $C > 0$ such that
\be
 \| u \, - \, u_{N} \|_{\omega^{-1}} \ \le \   \frac{1}{ \lambda_{N+1} } \, \| f \|_{\omega^{*}} \ \le \
 \ C \, (N + 2)^{- \alpha}    \, \| f \|_{\omega^{*}} . 
  \label{apertg2}
\ee 
\label{thmAP2}
\end{theorem} 
\textbf{Proof}: The proof follows in a similar manner to that for Theorem \ref{thmAP1}. \\
\mbox{ } \hfill \qed

\begin{corollary}
For $f(x) \in L^{2}_{\omega^{*}}(0 , 1)$ and $u_{N}(x)$ given by  \eqref{defuNs}, there exists $C > 0$ such that
\be
 \| u \, - \, u_{N} \|  \ \le \   \frac{1}{\lambda_{N+1} } \, \| f \|_{\omega^{*}} \ \le \
 \ C \, (N + 2)^{- \alpha }    \, \| f \|_{\omega^{*}} \, . \label{apertg3}
\ee 
\end{corollary} 
\label{corAPg1}
\textbf{Proof}: As $\omega(x) \, = \, x^{\beta} (1 - x)^{\alpha - \beta} \, < 1 $, for $0 < x < 1$, then
$\| u \, - \, u_{N} \| \ \le \ \| u \, - \, u_{N} \|_{\omega^{-1}}$. Hence 
the bound \eqref{apertg3} follows immediately from \eqref{apertg2} .
  \\
\mbox{ } \hfill \qed

\subsection{Numerical Examples}
\label{ssecspNE}
In this section we demonstrate the spectral type approximation methods discussed in Sections \ref{ssecsp1} and
\ref{ssecsp2} on Examples 1 and 2 presented in Section \ref{secFEMcvg}.

\begin{figure}[!ht]
\begin{minipage}{.46\linewidth}

\begin{center}
 \includegraphics[height=2.5in]{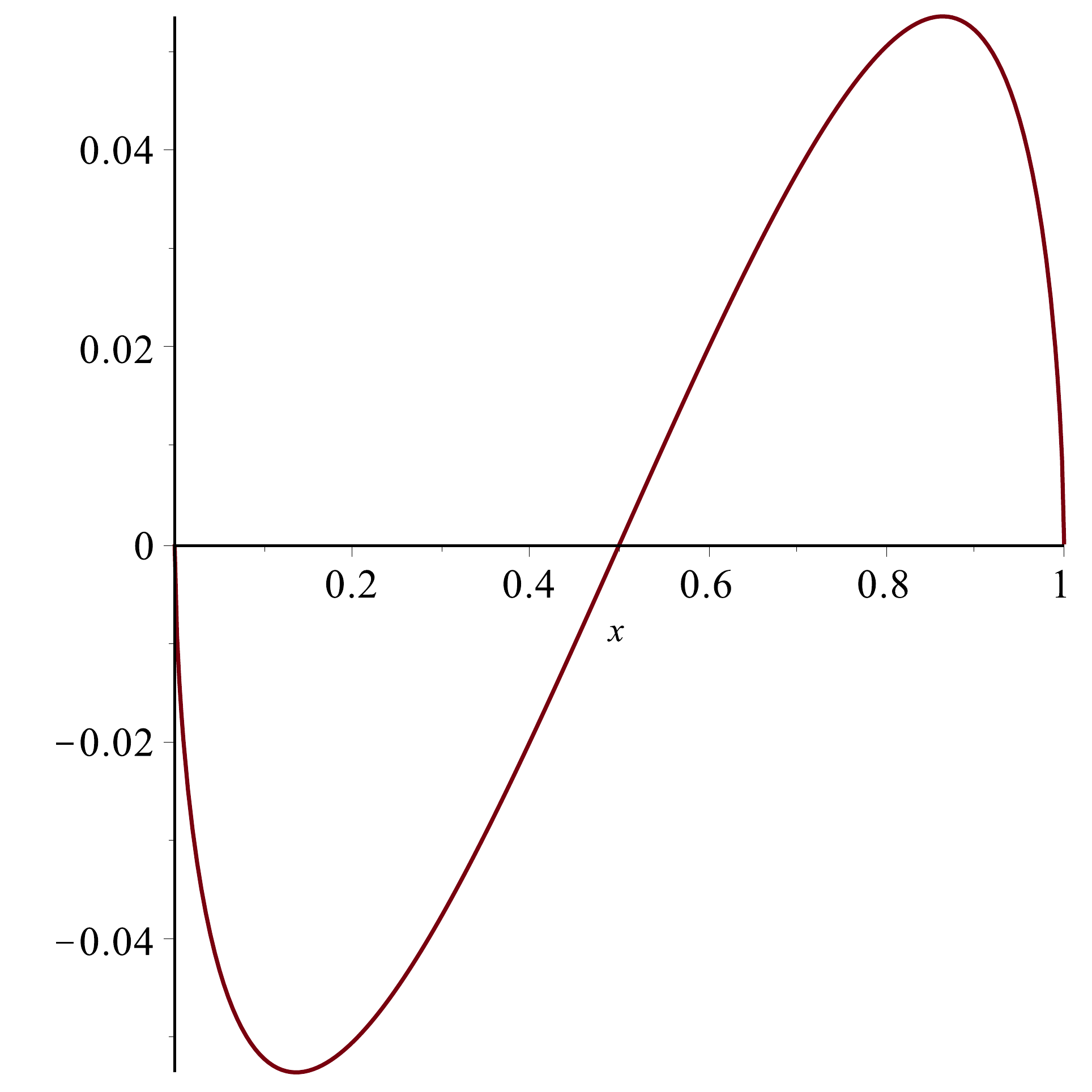}
   \caption{Solution of Example 1, $u(x)$ given in \eqref{ppo10}.}
   \label{spfig1}
\end{center}

\end{minipage} \hfill
\begin{minipage}{.46\linewidth}
 
\begin{center}
\includegraphics[height=2.5in]{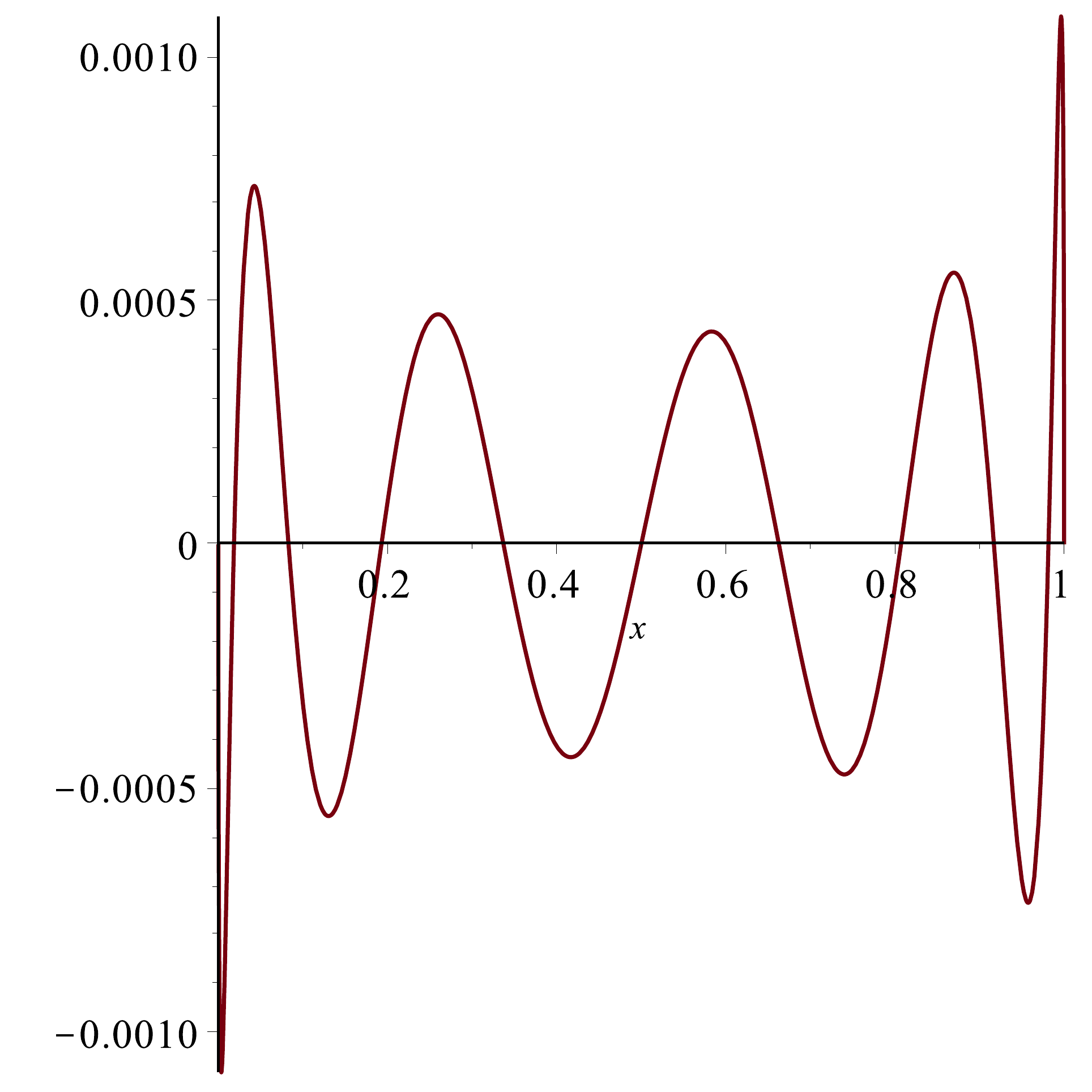}
   \caption{Plot of $u(x) \, - \, u_{8}(x)$ for Example~1.}
   \label{spfig2}
\end{center}
\end{minipage} 
\end{figure}

\begin{figure}[!ht]
\begin{minipage}{.46\linewidth}

\begin{center}
 \includegraphics[height=2.5in]{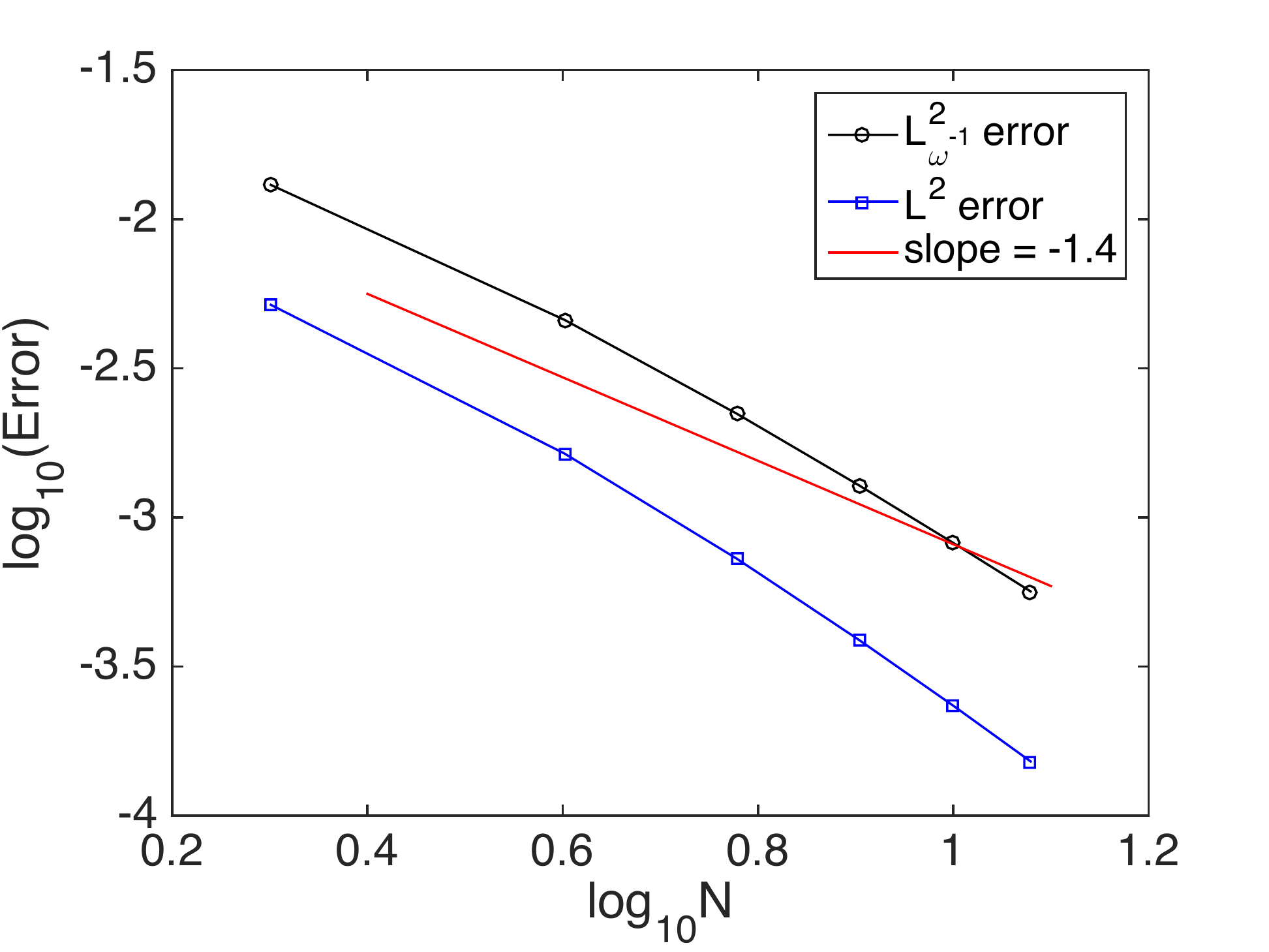}
   \caption{$L^{2}_{\omega}$ and $L^{2}$ errors for Example 1.}
   \label{spfig3}
\end{center}

\end{minipage} \hfill
\begin{minipage}{.46\linewidth}
 
\begin{center}
 \includegraphics[height=2.5in]{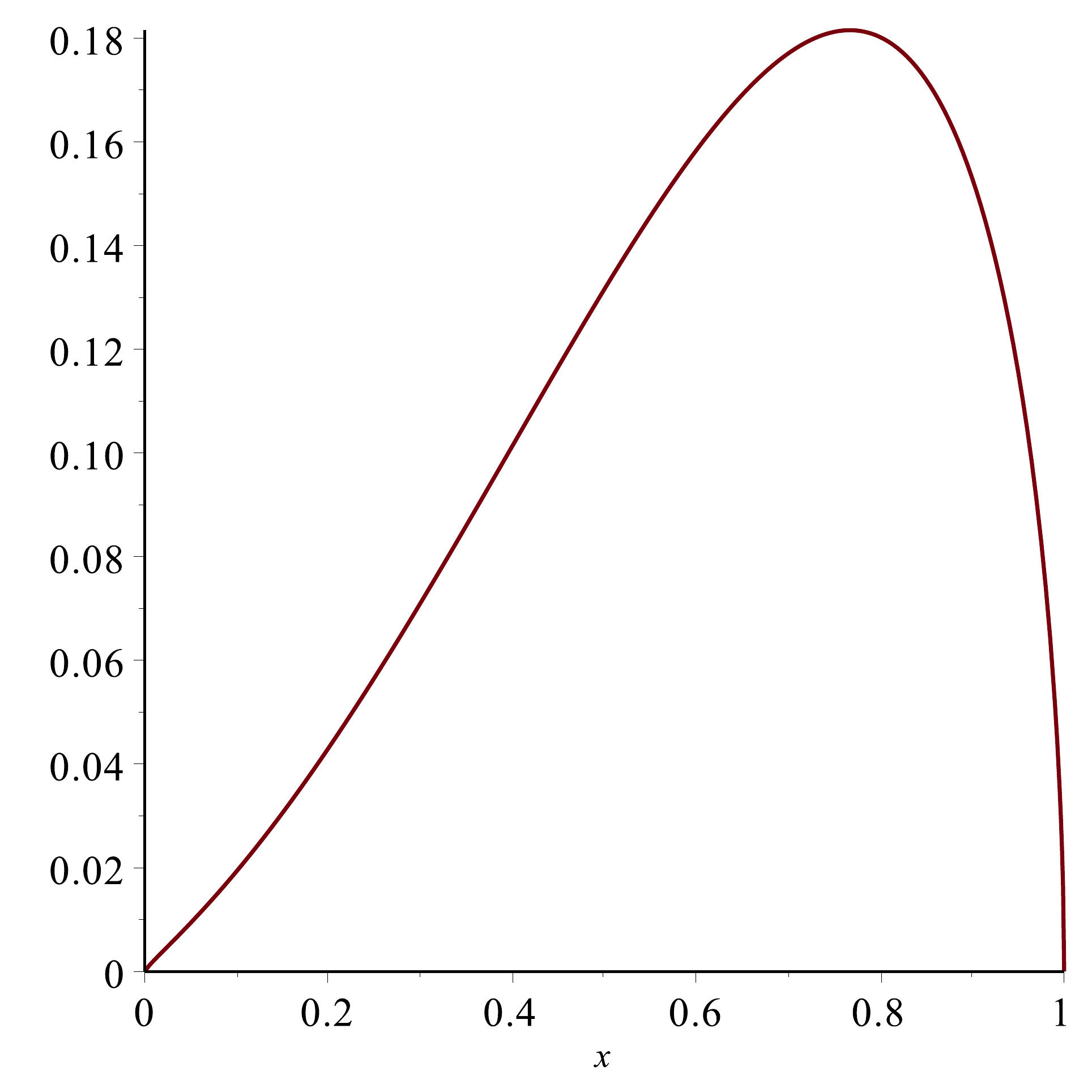}
   \caption{Solution of Example 2, $u(x)$ given in \eqref{ppo12}.}
   \label{spfig4}
\end{center}
\end{minipage} 
\end{figure}

\underline{Example 1. cont.} \\
For this example $\alpha = 1.4$ and $r = 1/2$. Hence we have (from Section \ref{ssecsp1}) that
$\omega(x) \, = \, x^{\alpha/2} (1 - x)^{\alpha/2} \, = \, x^{0.7} (1 - x)^{0.7}$, and from \eqref{deffi} and
\eqref{defuN}
\[
    u_{N}(x) \ = \ x^{0.7} (1 - x)^{0.7} \, \sum_{j = 0}^{N} \frac{f_{j}}{\lambda_{j} \, \| G_{j} \|_{w}^{2}} \, G_{j}(x)  \, .
\]
Presented in Figure \ref{spfig1} is a plot of the true solution given in \eqref{ppo10}. Figure \ref{spfig2} contains a
plot of the error, $u(x) \, - \, u_{8}(x)$, which exhibits a Gibbs type phenomena at the endpoints. Presented in 
Figure \ref{spfig3} is a plot of the $L^{2}_{\omega}$ and $L^{2}$ errors for the approximations. The 
convergence of the approximations  is consistent with the theoretical results given in \eqref{apert1} and
\eqref{apert3}. 

\begin{figure}[!ht]
\begin{minipage}{.46\linewidth}
 
\begin{center}
\includegraphics[height=2.5in]{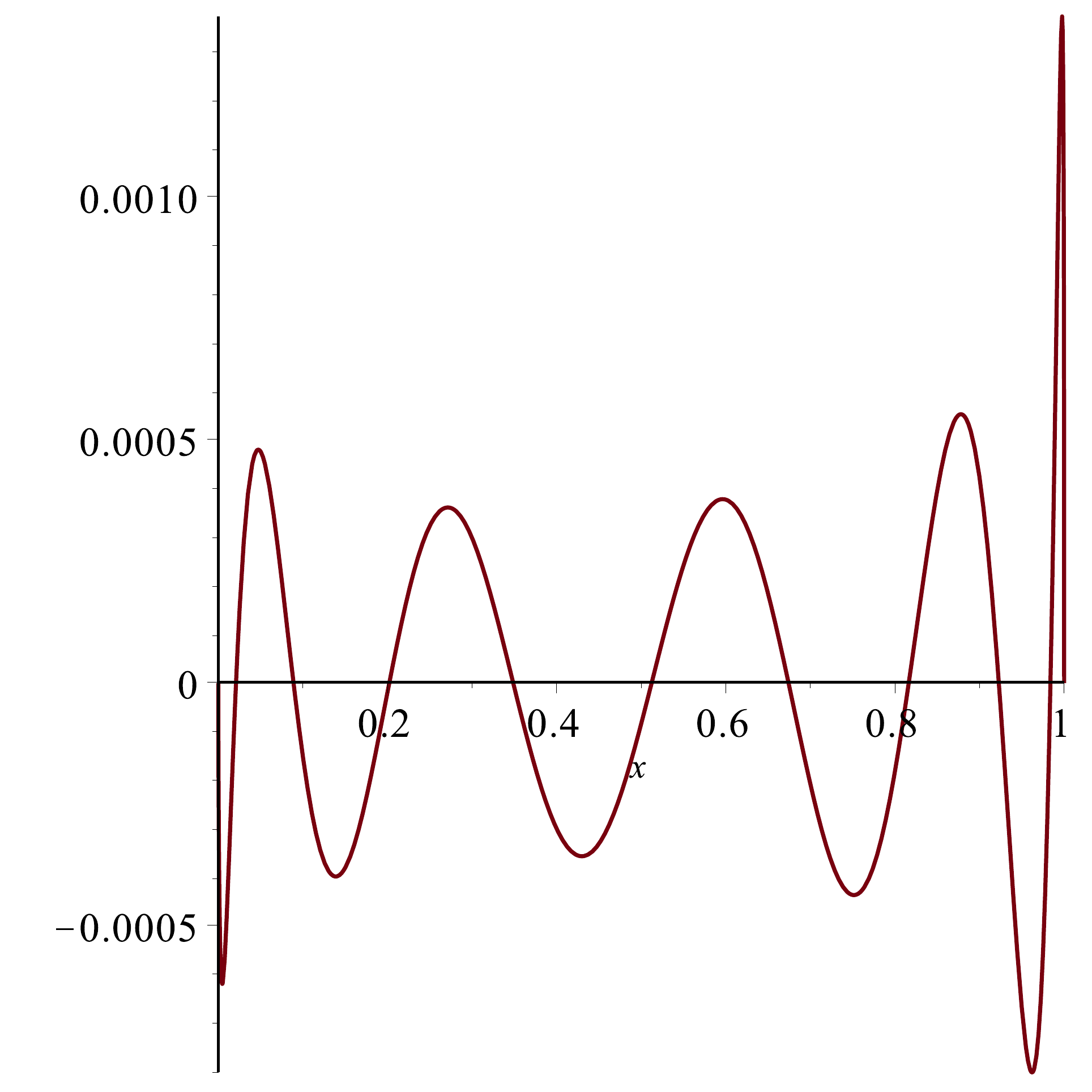}
   \caption{Plot of $u(x) \, - \, u_{8}(x)$ for Example~2.}
   \label{spfig5}
\end{center}
\end{minipage}  \hfill
%
\begin{minipage}{.46\linewidth}

\begin{center}
 \includegraphics[height=2.5in]{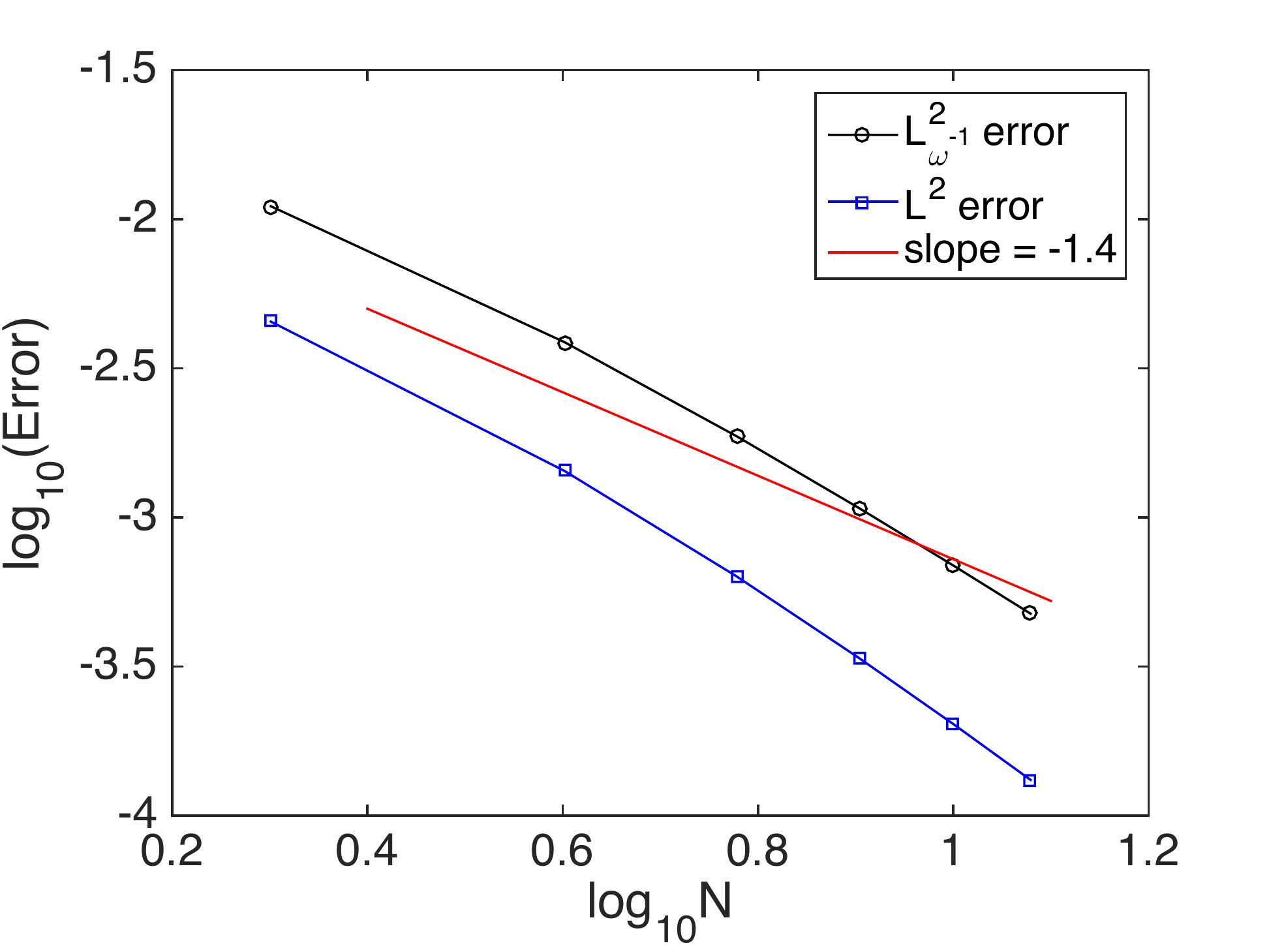}
   \caption{$L^{2}_{\omega}$ and $L^{2}$ errors for Example 2.}
   \label{spfig6}
\end{center}

\end{minipage} \hfill
\end{figure}

\underline{Example 2. cont.} \\
For this example $\alpha = 1.4$, $r = 0.3149$, $p = -0.15$ and $q = -0.45$. For these values the corresponding
value for $\beta = 0.85$ (see \eqref{spm11}). From Section \ref{ssecsp2}, \eqref{wwstar}, 
$\omega(x)  =  x^{\beta} (1 - x)^{\alpha - \beta} = x^{0.85} (1 - x)^{0.55}$, and from \eqref{deffis} and \eqref{defuNs}
\[
    u_{N}(x) \ = \ x^{0.85} (1 - x)^{0.55} \, \sum_{j = 0}^{N} \frac{f_{j}^{*}}{\lambda_{j} \, 
        \| G_{j}^{*} \|_{w^{*}}^{2}} \, G_{j}(x)  \, .
\]
Presented in Figure \ref{spfig4} is a plot of the true solution given in \eqref{ppo12}. Figure \ref{spfig5} contains a
plot of the error, $u(x) \, - \, u_{8}(x)$. Presented in 
Figure \ref{spfig6} is a plot of the $L^{2}_{\omega}$ and $L^{2}$ errors for the approximations. The 
convergence of the approximations  is consistent with the theoretical results given in \eqref{apertg2} and
\eqref{apertg3}. 

\underline{Acknowledgement}: The authors gratefully acknowledge helpful discussions with Professor
Jeff Geronimo on the properties of the hypergeometric functions.

\setcounter{equation}{0}
\setcounter{figure}{0}
\setcounter{table}{0}
\setcounter{theorem}{0}
\setcounter{lemma}{0}
\setcounter{corollary}{0}
\newpage
\begin{appendix}
\Large{\textbf{Appendix}}
\begin{figure}[!ht]  
\begin{center}
 \includegraphics[height=7.5in]{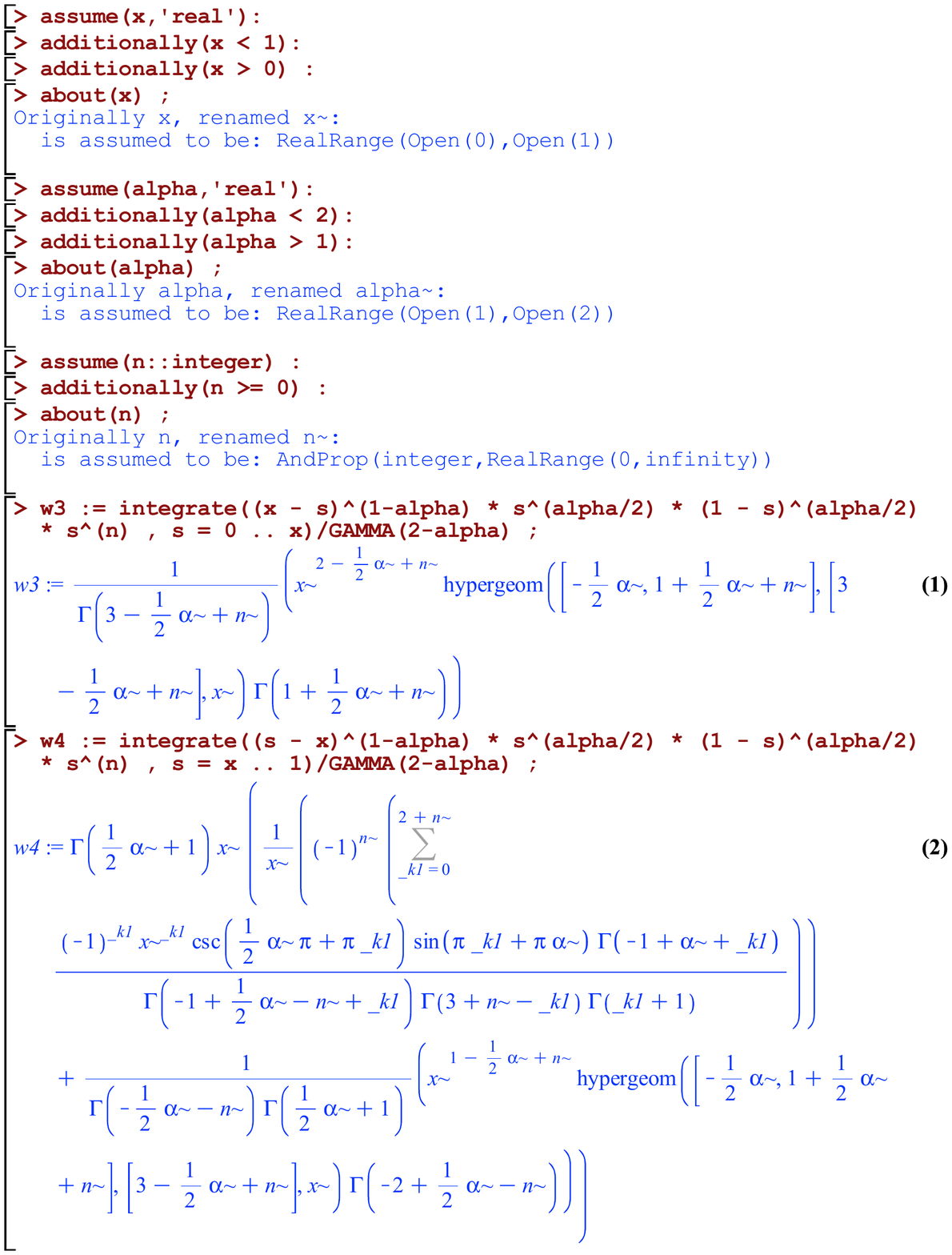}
   \caption{Maple computation for Lemma \ref{lmaSpec1}}
   \label{mapleSPL1}
\end{center}
\end{figure}

\begin{figure}[!ht]  
\begin{center}
 \includegraphics[height=7.5in]{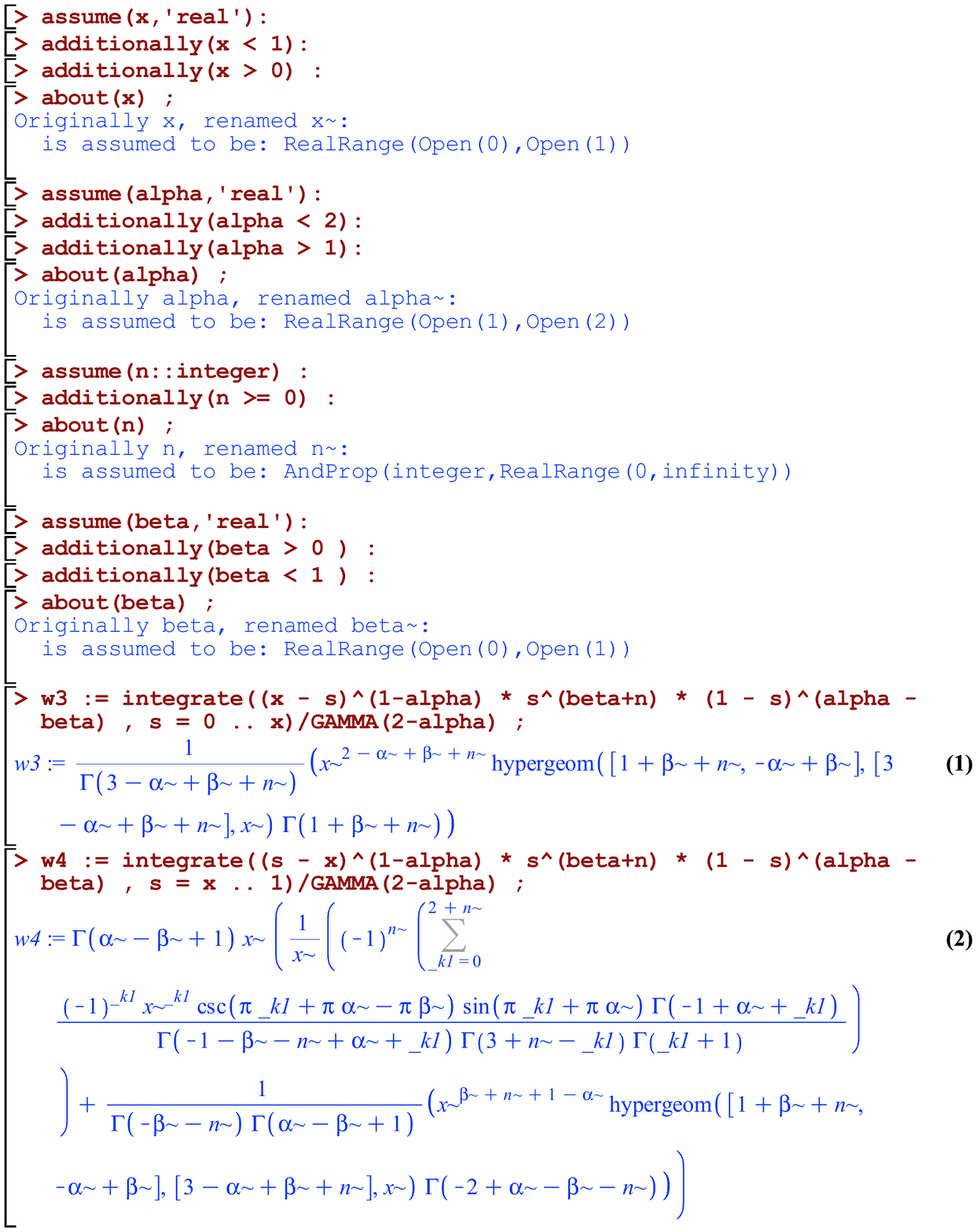}
   \caption{Maple computation for Lemma \ref{lmaSpec2}}
   \label{mapleSPL2}
\end{center}
\end{figure}

\end{appendix}

\end{document}